\documentclass[12pt]{article}
\usepackage{WideMargins2}
\def\Gal{{{\rm Gal}}}
\def\PGL{{{\rm PGL}}}

\def\SF{{{\rm SF}}}
\def\GCD{{{\rm GCD}}}
\def\Norm{{{\rm N}}}
\def\P{{{\mathbb P}}}
\def\F{{{\mathbb F}}}
\def\KK{{{\mathbb K}}}
\def\LL{{{\mathbb L}}}
\def\cj#1{{{\overline{#1}}}}

\def\deg{{{\rm deg}}}
\def\ang#1{{{\langle #1 \rangle}}} 
\def\textmatrix#1#2#3#4{{ \left({#1 \atop #3}{#2\atop #4}\right)}}
\let\a=\alpha

\let\g=\gamma

\let\l=\lambda

\let\x=\times
\let\s=\sigma
\let\z=\zeta
\def\ie{{ \it i.e.}}
\def\calC{{{\cal C}}}
\def\calD{{{\cal D}}}
\def\calY{{{\cal Y}}}

\title{\Large\bfseries  A New Identity of Dickson Polynomials }
\author{Antonia W.~Bluher}
\date{\today}
\begin {document}
\fancytitle

\begin{abstract}
A new polynomial identity is found for Dickson polynomials in characteristic~2.
The identity is used to prove that the two polynomials $x^{q+1}+x+1/a$ and
$C(x)+a$ have the same splitting fields over~$\F$,
where $\F$ is a field of characteristic~2, $0\ne a \in \F$,
$q = 2^n>2$, and $C(x)$ is a 
M\"uller--Cohen--Matthews polynomial of degree $(q^2-q)/2$.
In addition, a new proof  is obtained of the known result 
that $C(x)$ induces a permutation on $\F_{2^m}$ if and only if $2m$ and 
$n$ are relatively prime.
\end{abstract}

\section{Introduction} \label{sec:Intro}

We find a new polynomial identity in characteristic~2: 
\begin{equation}
\prod_{w\in\F_q^\x} (D_{q+1}(wX)-Y)
= X^{q^2-1} + \left(\sum_{i=1}^{n} Y^{2^{n}-2^i}\right) X^{q-1} + Y^{q-1},
\label{DicksonEq}
\end{equation}
where $q = 2^n$ and $D_k$ is a Dickson polynomial, 
defined by $D_k(u+u^{-1})=u^k + u^{-k}$.
Using this identity and results of \cite{qplus1}, we prove that if 
$\F$ is a field of characteristic~2
and $a$ is a nonzero element of $\F$, then for $q=2^n>2$, the two polynomials
$x^{q+1}+x+1/a$ and $C(x)+a$ have the same splitting field over~$\F$,
where $C(x) = x (\sum_{i=0}^{n-1} x^{2^i-1})^{q+1}$ is a 
M\"uller--Cohen--Matthews polynomial of degree $(q^2-q)/2$.
Explicit formulas are found relating the roots of the two polynomials,
and the Galois action on the roots of both polynomials is described.
As a result, when $\F$ is finite, related factorizations of the
two polynomials can be explained.
We also found a new proof that $C(x)$ induces a permutation on $\F_{2^m}$
if and only if $(2m,n)=1$.
(See \cite{CM} for the original proof. A polynomial over a finite field
$\F_\ell$ that induces a permutation
on infinitely many algebraic extension fields $\F_{\ell^k}$ is said to be 
{\it exceptional}. 

A first draft of this article was written
in the 2001--2004 time frame, but it was left and forgotten
for over a decade.  The project was resumed and completed in~2016, 
with the following
improvements. A hypothesis that the field $\F$ must be 
perfect\footnote{Since function fields are not perfect, it is important to 
remove the hypothesis that $\F$ be perfect.}
	was removed,
a new proof of exceptionality of $C(x)$ was added, a simpler formula 
for the roots of $C(x)$ was found, and
a simpler description of dihedral subgroups of $\PGL_2(q)$ was obtained.
References were updated (\cite{GRZ,GZ}) to reflect advances in 
the understanding of
exceptional polynomials that occurred in the intervening decade. 
The current version contains simplified proofs, particularly in
Section~\ref{qplus1Sec}.

A few remarks are in order.  First, equality of the splitting fields
of the two polynomials $x^{q+1}+x+1/a$ and $C(x)+a$ 
can be derived from work of Zieve \cite{Z} and
Lenstra and Zieve \cite{LZ}, at least in the case where $a$ is transcendental.
Many calculations in this paper could
perhaps be done more expediently with their methods, which utilize
group theory.  However, the author was unaware of these methods at the time that
she carried out her work, and as a result she used different techniques
and was motivated by a different set of questions.  We hope that this
new perspective will complement the existing literature.

The polynomial identity involving Dickson polynomials in characteristic~2 
is new.  It seems to apply only to characteristic~2.  
Bob Guralnick points out that
the Dickson polynomials are ramified at the prime~2, thus it is not
surprising to find formulas that are special to characteristic~2.

Our results seem related to but different from results of
Abhyankar, Cohen, and Zieve \cite{ACZ}.   Both our paper and theirs 
give a factorization of $x^{q^2-1}-a(y)x^{q-1}-b(y)$ in terms of
Dickson polynomials and use it to deduce information about the
Galois groups of certain polynomials. However,
the functions $a(y)$ and $b(y)$ differ, and so do the 
Galois groups that are involved.  
Their identity generalizes to all characteristics, whereas
ours applies only to characteristic~2.
A precise statement of the identity in \cite{ACZ} is given in
the remark preceding Lemma~\ref{DicksonIdLemma}.  It would be interesting
to understand more fully how the two polynomial identities are related.

Finally, we mention that in
one case, our work fits nicely with results of Dummit \cite{Dummit}
on solvable quintics. Namely,
$x^{q+1}+x+1/a$ is a quintic when $q=4$.
Dummit notes that a quintic (over any field)
is solvable if and only if its Galois group is
contained in a group that is conjugate to $F_{20}$, where
$F_{20}\subset S_5$ is generated by the permutations $(1 2 3 4 5)$
and $(2354)$. It turns out that an invariant $\theta$ for $F_{20}$ is
given by:
\begin{eqnarray*}
\theta &=& x_1^2 x_2 x_5 + x_1^2 x_3 x_4 + x_2^2 x_1 x_3 + x_2^2 x_4 x_5
+ x_3^2 x_1 x_5 \\
&& + x_3^2 x_2 x_4 + x_4^2 x_1 x_2 + x_4^2 x_3 x_5 + x_5^2 x_1 x_4 + 
x_5^2 x_2 x_3,
\end{eqnarray*}
where $x_i$ are the roots of the quintic.
Let $\g_1,\ldots,\g_6$ be coset representatives for $S_5/F_{20}$.
Then $\prod_{i=1}^6 (x - \g_i(\theta))$ is a sextic. In the case
of the polynomial $x^5+x+1/a$ in characteristic~2, the sextic is 
equal to $x^6+a^{-4}x+a^{-4}$. Setting $y=1/x$, 
we see that this is equivalent (has the same splitting field) to $y^6+y^5+a^4$.
If we set $y=z^4$ and then take a fourth root, this reduces to 
$z^5(z+1) + a$. The substitution $z \mapsto z+1$ brings
it to the form $C(z)+a$.
The above expression for $\theta$ shows explicitly in this case how
the roots of $C(x)+a$ are related to the roots of the quintic.

The relation between the roots causes a relation between the factorizations
of the two polynomials.  
Write $f \sim [n_1,n_2,...,n_t]$ if $f$ factors into irreducibles 
of degrees $n_1, n_2, \ldots , n_t$.
When $q=4$ and $\F=\F_{2^k}$, we will prove the following in Section~\ref{galoisSec}.

\begin{proposition} \label{quinticSextic}
\noindent For $k$ even, $x^5+x+1/a$ has one of these factorization types:
$[1,1,1,1,1]$, $[1,1,3]$, $[1,2,2]$, or $[5]$.  We have
\begin{eqnarray*}
x^5 + x + 1/a \sim [1,1,1,1,1] &\iff& x(1+x)^5 + a \sim [1,1,1,1,1,1]        \\
x^5 + x + 1/a \sim [1,1,3] &\iff& x(1+x)^5 + a \sim [3,3]        \\
x^5 + x + 1/a \sim [1,2,2] &\iff& x(1+x)^5 + a \sim [1,1,2,2]        \\
x^5 + x + 1/a \sim [5] &\iff& x(1+x)^5 + a \sim [1,5].       
\end{eqnarray*}
For $k$ odd, $x^5+x+1/a$ has factorization type $[1,1,1,2]$, $[1,4]$, or $[2,3]$, and 
\begin{eqnarray*}
x^5 + x + 1/a \sim [1,1,1,2] &\iff& x(1+x)^5 + a \sim [2,2,2]       \\
x^5 + x + 1/a \sim [1,4] &\iff& x(1+x)^5 + a \sim [1,1,4]       \\
x^5 + x + 1/a \sim [2,3] &\iff& x(1+x)^5 + a \sim [6].       
\end{eqnarray*}
\end{proposition}

The author wishes to thank John Dillon and Mike Zieve for
their encouragement and for stimulating discussions.
As a newcomer to the study of Galois groups and exceptional polynomials,
the author found their expertise to be invaluable.

\medskip
\noindent {\bf Notation.}
If $\F$ is a field, then $\F^\times$ denotes its group of units
and $\cj\F$ denotes its algebraic closure.
If $f\in \F[x]$, the splitting field of $f$ over~$\F$ is 
written $\SF(f;\F)$; this is the
subfield of $\cj \F$ that is generated by $\F$ and by all the roots of $f$.
The Galois group of $\SF(f;\F)$ over $\F$ is denoted by $\Gal(f;\F)$; it is 
the group of automorphisms of the field $\SF(f;\F)$ that fix the subfield~$\F$.
If $\ell$ is a prime power, then $\F_\ell$ denotes the (unique) field
with $\ell$ elements. 
We sometimes abbreviate ``if and only if'' by ``iff'' or the symbol $\iff$.

For $k\ge 1$, $D_k(x)$ denotes the $k$th Dickson
polynomial, which is the unique monic polynomial of degree $k$ such that
$D_k(x+1/x)=x^k + 1/x^k$.  
Because the expression $x+1/x$ arises frequently, we introduce
the special notation:
$$\langle x \rangle = x + 1/x.$$
Then the defining property of the Dickson polynomial may be written as
$$D_k(\langle x \rangle) = \langle x^k \rangle.$$
Note that
\begin{eqnarray}
\ang {1/x} &=& \ang x \nonumber \\
	\label{xyangProperty}
\langle x \rangle  \langle y \rangle &=& \langle x y \rangle
+ \langle x/y \rangle  \\
\label{angProperty}
\ang x \ang{y/z} &=& \ang{xy}\ang z + \ang{xz} \ang y \qquad\text{in char.~2} 
\\
\langle x^{p^i} \rangle &=& \langle x \rangle^{p^i}\qquad\text{ in char.~$p$.}
\end{eqnarray}

The following notation is specific to this article:
$$q = 2^n, \qquad{\rm where\ }n > 1.$$
$$T(x) = \sum_{i=0}^{n-1} x^{2^i-1}$$
%T_{rev}(x)=\sum_{i=0}^{n-1}x^{(q/2)-2^i}$$
$$C(x) = x\cdot T(x)^{q+1}\quad\text{(a M\"uller--Cohen--Matthews polynomial)}$$
$$\F\ \text{is a field of char.~2 and $a$ is a fixed nonzero element of $\F$}$$
$$\text{If $k\ge1$ is odd, then 
$\mu_{k} = \{\zeta\in\cj \F_2^\x : \zeta^{k}=1\}$}$$
$$ \KK=\SF(C(x)+a;\F),\qquad \LL = \SF(x^{q+1}+ x + 1/a;\F)$$
$$\F_{q,1} = \{\,d \in \F_q : \Tr_{\F_q/\F_2}(d) = 1 \},\qquad 
\F_{q,0} = \{\,d \in \F_q : \Tr_{\F_q/\F_2}(d) = 0 \}.$$

\noindent In Sections~\ref{sec:Identity} and ~\ref{qplus1Sec} only, 
we allow $q=p^n$ and $\F$ has char.~$p$, where $p$ is any prime and $n\ge 1$.
We make frequent use of elements $\zeta\in\mu_{q+1}\subset \mu_{q^2-1} =
\F_{q^2}^\x$.  Note that
$\Norm_{\F_{q^2}/\F_q}(\zeta)=\zeta\zeta^q=1$,
and $\Tr_{\F_{q^2}/\F_q}(\zeta)=\zeta+\zeta^q=\zeta+\zeta^{-1}=\ang\zeta$.

%\medskip
%\noindent{\bf Remark.}  
%If one assumes that 
%$\F$ is perfect, then replacing all the coefficients of $f$ by their $2^i$th powers
%does not affect the splitting field. For example, 
%the substitution $r \mapsto r/a$ transforms roots of
%$x^{q+1} + x + 1/a$ to roots of $x^{q+1} + a^q x + a^q$, and if $\F$ is
%perfect then this has the same splitting field 
%as $x^{q+1} + a x + a$.  
%In this article, we do not assume that $\F$ is perfect, but we are still able
%to show in Proposition~\ref{bpLemma} that the polynomials
%$x^{q+1} + x + 1/a$ and $x^{q+1} + ax + a$ have the same splitting field.
%The latter polynomial turns out to be
%the most convenient for the purpose of proving that these polynomials have the same
%splitting field as $C(x)+a$. 

%\smallskip
%The following lemma will be needed in Sections~\ref{sec:RootsOfC}
%and~\ref{dihedralSec}.

\begin{lemma}  {\rm (Dillon and Dobbertin, \cite[pp.~154--155]{DD})} \label{Fq1Lemma}
Let $q=2^n$. Then
$$\F_{q,1}=\{1/\ang\zeta : \zeta \in \mu_{q+1}, \zeta\ne 1 \},\qquad
\F_{q,0}=\{0\}\cup \{1/\ang w : w \in \mu_{q-1}, w\ne 1 \}.$$ 
If $0\ne a \in \F_q$, then 
\begin{equation*} \text{$x^2+ax+1$ is reducible over $\F_q$
iff $1/a\in\F_{q,0}$ iff $a=\ang w$ with $w\in \F_q^\x\setminus\{1\}$;}\
\end{equation*}
\begin{equation*} \text{$x^2+ax+1$ is irreducible over $\F_q$ iff 
$1/a\in\F_{q,1}$ iff $a=\ang \zeta$ with $\zeta\in \mu_{q+1}\setminus\{1\}$.}
\end{equation*}
\end{lemma}
\begin{proof}
Write $x^2+ax+1=(x-r)(x-1/r)$, where $r\in\F_{q^2}^\x$. Then $a=\ang{r}$
and $(r/a)^2+(r/a)=1/a^2$. 
Thus, $x^2+ax+1$ is reducible iff $r\in\F_q^\x$ iff $1/a^2= b^2+b$
has a solution $b=r/a\in\F_q$ iff $\Tr(1/a)=0$.

If $x^2+ax+1$ is irreducible, then its roots form a conjugate pair, 
so $1/r=r^q$.
Thus, $a=\ang r$ with $r\in\mu_{q+1}\setminus \{1\}$.
Since the solutions to $1/a^2 = b^2 + b$ are irrational, $\Tr(1/a)=1$.
%
%\noindent
	%({\it ii}) Let $f(\rho)=\ang{\z/\rho}/(\ang\z\ang\rho)$.
	%Since $\ang\lambda \in \F_q$ for all $\lambda \in \mu_{q+1}$,
	%and $\ang\lambda=0$ iff $\lambda=1$,  $f$
	%maps $\mu_{q+1}\setminus \{1\}$ into $\F_q$.
	%The two sets $\mu_{q+1}\setminus\{1\}$ 
%and $\F_q$ have the same cardinality, so it suffices to prove $f$ is injective.
	%Suppose $\rho,\rho'\in \mu_{q+1}\setminus \{1\}$ and
	%$f(\rho)=f(\rho')$. 
%We will show $\rho=\rho'$. The equality 
%implies $\ang{\z/\rho}\ang{\rho'} = \ang{\z/\rho'}\ang\rho$. 
	%Applying (\ref{angProperty}) with $x=\z$, $y=1/\rho$, $z=1/\rho'$ shows
	%$\ang\z\ang{\rho'/\rho}=0$. Since $\ang\z\ne0$, it must be
	%that $\ang{\rho'/\rho}=0$. Then $\rho'=\rho$, as required.

The formulas for $\F_{q,0}$ and $\F_{q,1}$ follow from these observations.
\end{proof}
The article is organized as follows.
Section~\ref{sec:Identity} proves the Dickson polynomial identity.  Section~\ref{qplus1Sec} concerns $x^{q+1}+ax+a$.
It reviews \cite{qplus1} and proves a few additional results, including
that SF$(x^{q+1}+ax+a) = {\rm SF}(x^{q+1}+x+1/a)=\LL$.
Sections~\ref{sec:RootsOfC} and~\ref{sec:Equality} prove that
$\KK = \LL$ by explicitly writing the roots of $C(x)+a$
as a rational function of the roots of $x^{q+1}+ax+a$, and explicitly writing the roots
of $x^{q+1}+ax+a$ as a rational function of the roots of $C(x)+a$. 
%Also the following polynomial identity is derived:
%$$\prod_{c\in\F_q^\x,\ j \in \F_{q,1}} (c y^2+y+j/c) = 1 + (y^q+y)^{q-1}.$$
Section~\ref{galoisSec} considers the Galois group and shows how the factorizations
of $x^{q+1}+x+1/a$ and $C(x)+a$ are related. For example, we prove the related
factorizations between $x^5 + x + 1/a$ and $x(x+1)^5 + a$ that were asserted in Proposition~\ref{quinticSextic}.
Section~\ref{dihedralSec} investigates dihedral groups of order $2(q+1)$ and shows that such groups
fix a root of $C(x)+a$ in the geometric case. This is used in Section~\ref{sec:Exceptional} to give a new proof that
$C(x)$ is exceptional over $\F_2$ when $n$  is odd.

\section{An identity of Dickson polynomials} \label{sec:Identity}

The $k$th Dickson polynomial 
is the polynomial with integer coefficients such that the
formal identity holds, $D_k(u+1/u) = u^k + 1/u^k$.  To see that
such a polynomial exists, note that $u^k +v^k$ is a polynomial in
$u+v$ and $uv$ by the Theorem of Symmetric Functions, say 
$u^k+v^k=F_k(u+v,uv)$. By setting $v=1/u$, we find that
$D_k(x) = F_k(x,1)$. It is easy to see that $D_k(x)$ is monic and
has degree $k$.  A useful relation is
\begin{equation}
\label{DkDl}
D_k(x) D_\ell(x) = D_{k+\ell}(x) + D_{k-\ell}(x)
\end{equation}
when $k\ge\ell$, as can be seen from the identity $\ang{u^k}\ang{u^\ell} =
\ang{u^{k+\ell}} + \ang{u^{k-\ell}}$.

Since $D_k$ has integral coefficients, it can be considered over
any field $\F$, in any characteristic.  If the characteristic is~$p$,
then
\begin{equation}
\label{Dkp}
D_{kp^r}(x) = D_k(x)^{p^r}.
\end{equation}
The complete set of roots of $D_k(x) - c$ is easy to construct.
Namely, we find $u$ so that $c=\ang{u^k}$; then $\ang{u}$ will be a root,
and the other roots will be $\ang{\zeta u}$ for $\zeta\in\mu_k$.
To find $u$, first solve the quadratic $v+1/v = c$, then solve $v=u^k$.

We will need some well-known formulas for $D_{q-1}$ and $D_{q+1}$ in 
characteristic~2, where $q=2^n$. For the reader's convenience,
their proof is included below.

\begin{lemma}
\label{lem:Dicksonqm1}
If $q=p^n$, where $n\ge 1$, then  in characteristic $p$,
$$D_{q+1}(Y) = Y^{q+1} - D_{q-1}(Y).$$ 
If $q=2^n$, then  in characteristic~2,
$$D_{q-1}(Y) = \sum_{i=1}^n Y^{q-2^i+1}.$$
\end{lemma}

\begin{proof}
By~(\ref{DkDl}), $D_{q}(Y)D_1(Y)=D_{q+1}(Y)+D_{q-1}(Y)$.
By~(\ref{Dkp}), $D_q(Y)D_1(Y)=Y^{q+1}$. Thus, $D_{q+1}(Y)=Y^{q+1}-D_{q-1}(Y)$.

Now apply~(\ref{DkDl}) with $k=q-1$ and $\ell=1$.
We find that $YD_{q-1}(Y) = Y^q + D_{q-2}(Y)$. If $p=2$, this becomes
$$Y^q = Y D_{q-1}(Y) + D_{(q/2)-1}(Y)^2.$$
Let $f_m = D_{2^m-1}(Y)$ and $g_m =\sum_{i=1}^m Y^{2^m-2^i+1}$. 
Then $f_1=g_1=Y$, $Yf_{m}=Y^{2^m} + f_{m-1}^2$ for $m\ge2$, and
$Yg_{m}=Y^{2^m} + g_{m-1}^2$ for $m\ge2$.  Thus, $f_n=g_n$.
\end{proof}

\begin{lemma} \label{distinctLemma}  Let $q=2^n\ge2$, 
$M$ a field of characteristic~2 that strictly contains
$\F_{q^2}$, and $u\in M \setminus \F_{q^2}$. Then the values 
$\{\,\ang{\zeta u}/w :  \zeta\in\mu_{q+1} {\rm\ and\ } w\in\F_q^\x\,\}$ 
are distinct.
\end{lemma}

\begin{proof}   Let $\zeta,\lambda \in \mu_{q+1}$ and $w,w' \in \F_q^\x$, and suppose that $\ang{\zeta u}/w = \ang{\lambda u}/w'$.
Then $w/w'=\ang{\zeta u}/\ang{\lambda u}$, and so
$$\frac{\zeta w}{\lambda w'} = \frac{\zeta \ang{\zeta u}}{\lambda \ang{\lambda u}} = \frac{\zeta^2 u^2 + 1}{\lambda^2 u^2 + 1}.$$
If $\zeta \ne \lambda$, then we can solve for $u^2$ in terms of $\zeta,\lambda,w,w'$, but this contradicts the hypothesis
that  $u \not \in \F_{q^2}$.
Thus, $\zeta = \lambda$, and consequently $w=w'$ also. We have shown that $\ang{\zeta u}/w=\ang{\lambda u}/w'$ implies $\zeta=\lambda$ and
$w=w'$, so these values are distinct, as claimed.
\end{proof}

\begin{theorem} 
\label{DicksonIdentity}
Let $q=2^n$, where $n\ge 1$. In the polynomial ring $\F_q[X,Y]$, 
	the identity (\ref{DicksonEq}) holds.
%\begin{equation}
%\label{DicksonEq}
%\prod_{w\in\F_q^\x} (D_{q+1}(wX)-Y)
%= X^{q^2-1} + \left(\sum_{i=1}^{n} Y^{2^{n}-2^i}\right) X^{q-1} + Y^{q-1}.
%\end{equation}
\end{theorem}

\begin{proof}
Let $U$ be transcendental over $\F_2$ and $$Y=\ang{U^{q+1}}.$$ 
Then $Y$ is also transcendental, and $\F_2(Y)\subset\F_2(U)$.  
Let $L(X)$ and $R(X)$ denote the left-hand and right-hand sides
of (\ref{DicksonEq}) respectively, considered as elements of $\cj\F_2(U)[X]$. 
Both are monic polynomials in $X$ of degree
$q^2-1$, and so $\deg_X(L-R) < q^2-1$.
Thus, to prove $L-R$ is identically zero it suffices to find 
$q^2-1$ distinct roots in $\cj\F_2(U)$.  We claim that these roots are 
\begin{equation}\{\,\ang{\zeta U}/w :  \zeta\in\mu_{q+1} {\rm\ and\ } w\in\F_q^\x\,\}.
\label{identityRoots}
\end{equation}
(These are distinct by Lemma~\ref{distinctLemma}.)
In fact we will show that $L$ and $R$ each vanish at these values.
Let $x$ denote one of these values: $$x=\ang{\zeta U}/w.$$ 
First, $L(x)=0$ because
$$D_{q+1}(wx)-Y= D_{q+1}(\ang{\zeta U}) -Y
= \ang{(\zeta U)^{q+1}} - \ang{U^{q+1}} = 0.$$
Next we show $R(x)=0$. Set $V=\zeta U$; then 
\begin{equation}
x = \ang{V}/w \qquad{\rm and } \qquad Y = \ang{V^{q+1}}.
\label{wxY}
\end{equation}
Note that $\ang{V}Y=\ang{V}\ang{V^{q+1}}$ is nonzero, since $V$ is 
transcendental.  Thus, it will suffice to show that $\ang V Y\,R(x)=0$.
By Lemma~\ref{lem:Dicksonqm1}, 
\begin{eqnarray*}
\ang V Y\,R(x) &=& \ang V Y (\ang V/w)^{q^2-1} 
+\ang VY\left(\sum_{i=1}^n Y^{q-2^i}\right)
(\ang V/w)^{q-1} + \ang V Y^q \\
&=& Y{\ang V}^{q^2} + \left(\sum_{i=1}^n Y^{2^n-2^i+1}\right) {\ang V}^q 
+ \ang V Y^q  \\
&=& \ang{V^{q+1}} \ang{V^{q^2}} + D_{q-1}(Y) \ang{V^q} + \ang V \ang{V^{q(q+1)}}.
\end{eqnarray*} 
Now $D_{q-1}(Y)=D_{q-1}(\ang{V^{q+1}}) = \ang{(V^{q+1})^{q-1} } =
\ang{V^{q^2-1}}.$  Using this observation and (\ref{xyangProperty}), 
we obtain
\begin{eqnarray*}
\ang V Y\,R(x) 
&=& \ang{V^{q+1}} \ang{V^{q^2}} + \ang{V^{q^2-1}} \ang{V^q} + \ang V \ang{V^{q(q+1)}} \\
&=& \ang{V^{q^2+q+1}} + \ang{V^{q^2-q-1}} + \ang{V^{q^2+q-1}} + \ang{V^{q^2-q-1}} + \ang{V^{q^2+q+1}} + \ang{V^{q^2+q-1}} \\ 
&=&  0.
\end{eqnarray*} 
Thus, $R(x)=0$, as claimed.
\end{proof}

\noindent{\bf Remark.}  Our
identity~(\ref{DicksonEq}) is tantalizingly
similar to an identity in Theorem~(1.1) from the
article by Abhyankar, Cohen, and Zieve \cite{ACZ}.
Their identity is
$$X^{q^2-1} - E_q(Y,1) X^{q-1} + E_{q-1}(Y,1)
=(X^{2q-2} - Y X^{q-1} + 1)\left(
\prod_{w\in\F_q^\x} (D_{q-1}(X,w) - Y) \right),$$
where $q=p^n$, $D_n(X,a)$ is defined by $D_n(U_1+U_2,U_1U_2)=U_1^n+U_2^n$,
and $E_n(Y,a)$ is defined by 
$E_n(U_1+U_2,U_1U_2) = (U_1^{n+1}-U_2^{n+1})/(U_1-U_2)$.
Using the relations (2.20) and (2.9) of \cite{ACZ}, this identity can 
be rewritten when $p=2$ as:
$$X^{q^2-1} + (D_{q+1}(Y)/Y) X^{q-1} + Y^{q-1}
=(X^{2q-2} - Y X^{q-1} + 1)\left(
\prod_{w\in\F_q^\x} (D_{q-1}(wX) - Y) \right),$$
and using Lemma~2.1, our identity can be rewritten as
$$X^{q^2-1} + (D_{q-1}(Y)/Y) X^{q-1} + Y^{q-1}
= \prod_{w\in\F_q^\x} \left(D_{q+1}(wX) - Y \right).$$
Bob Guralnick points out that
the Dickson polynomials are ramified at the prime~2, so it is not
surprising to find formulas that are special to characteristic~2.

For future use, we record the following lemma.

\begin{lemma} \label{DicksonIdLemma}
	Let $q=2^n\ge 2$,
$y$ a nonzero element of a field $M$ of characteristic~2,
and 
$$f(x) = \prod_{w\in\F_q^\x}(D_{q+1}(wx)-y).$$
Let $u\in \cj M$ satisfy $u^{q+1} + 1/u^{q+1} = y$. 
The complete set of roots of $f$ in $\cj M$ is 
$$\{\,w \ang{\zeta u}  : w\in\F_q^\x, \zeta\in\mu_{q+1}\,\}$$
and these roots are distinct.
\end{lemma}

\begin{proof}   By (\ref{DicksonEq}),
	$f = x^{q^2-1} + \left(\sum_{i=1}^n y^{q-2^i}\right) x^{q-1} + y^{q-1}$.
%$f = x^{q^2-1} + T_{rev}(y^2) x^{q-1} + y^{q-1}$.
The roots of $f$ are distinct, because
$f-xf'=y^{q-1}\in M^\x$ implies $\GCD(f,f')=1$. 
Also the roots are nonzero,
since the constant term of $f$ is $y^{q-1}$.
Now $D_{q+1}(wx)-y$ vanishes at $w^{-1}\ang{\zeta u}$, for all $\zeta\in\mu_{q+1}$. 
We claim the values $\ang{\zeta u}$ are
distinct. If not, then $\ang{\zeta u} = \ang{\zeta' u}$ for distinct
$\zeta,\zeta'\in \mu_{q+1}$. One finds that $\zeta u^2 + 1/\z =
	\zeta' u^2 + 1/\zeta'$, therefore $u^2=(1/\zeta' + 1/\zeta)/(\zeta
+\zeta') = 1/(\zeta \zeta')$. Then
$y^2=\ang{u^{2(q+1)}}=\ang1=0$, contrary to the hypothesis
that $y$ is nonzero.  This establishes that the $q+1$ roots of $D_{q+1}(wx)-y$
given by $\{\ang{\zeta u}/w:\zeta\in\mu_{q+1}\}$ are distinct.
The roots of $D_{q+1}(wx)-y$ must be disjoint
from those of $D_{q+1}(w'x)-y$ when $w\ne w'$, since we already observed
that $f$ has no repeated roots.  Thus, the $q^2-1$ roots given in the statement
of the lemma are distinct, and since $\deg(f)=q^2-1$, we have found all the roots.
\end{proof}

Now we show how the Dickson polynomial identity leads to a relation between the polynomials
$C(x)+a$ and $x^{q+1} + a x + a$, where we recall from {\it Notation} that
$a$ is a nonzero element of a field $\F$ of char.~2, and $C$ is defined by
$$C(x) = x\cdot T(x)^{q+1},\qquad T(x)=\sum_{i=0}^{n-1} x^{2^i-1},
\qquad q = 2^n \ge 4.$$

\begin{lemma}  \label{lem:distinctRootsC}
$C(x)+a$ has distinct roots in $\cj \F$, all nonzero.
\end{lemma}

\begin{proof} A polynomial has distinct roots in the algebraic closure iff it is
relatively prime to its derivative.
Since $C=xT^{q+1}$ and $xT'=T+1$, we see that $C'=T^q(xT'+T)=T^q$.
Setting $G(x) = C(x)+a$, we have $G' = C' = T^q$,
	and so $G - x T G' = (C + a) - xT^{q+1} = a\in \F^\x$.
This proves that $\GCD(G,G')=1$, and so $C(x)+a$ has no repeated roots. 
Since $C(0)+a=a\ne0$, the roots are nonzero.
\end{proof}

For the remainder of this section, let $e$ be an arbitrary 
root of $C(x)+a$: 
$$ C(e) = a.$$

\begin{proposition}  
\label{eRelation}
Let $u\in\cj \F$ satisfy
$$u^{q+1} + 1/u^{q+1} = 1/e.$$
Then the complete set of roots of $x^{q+1} + a^2 x + a^2$ is
$$\left\{\,e T(e)^2 \ang{\zeta u }^{q-1}: \zeta\in\mu_{q+1}\,\right\}.$$
\end{proposition}

\begin{proof} Substitute $Y=1/e$ into the identity (\ref{DicksonEq}) 
and leave $X$ as an indeterminate.  We find:
$$\prod_{w\in\F_q^\times} (D_{q+1}(wX)-1/e) 
= X^{q^2-1} + \alpha X^{q-1} + \beta,$$
where
$$\alpha = \sum_{i=1}^n e^{2^i-q} = (e T(e))^2/e^q = e^{2-q} T(e)^2,\qquad \beta = e^{1-q}.$$
Let 
$$R = (\alpha/\beta) X^{q-1} = e T(e)^2 X^{q-1}.$$
The right side of the identity can be written as 
$(\beta/\alpha)^{q+1}(R^{q+1} + \alpha^{q+1}/\beta^q R + \alpha^{q+1}/\beta^q)$.
Now $\alpha^{q+1}/\beta^q = e^2 T(e)^{2(q+1)} = C(e)^2 = a^2$, therefore
$$\prod_{w\in\F_q^\times} (D_{q+1}(wX)-1/e) = (\beta/\alpha)^{q+1}
(R^{q+1}+a^2 R + a^2).$$
By Lemma~\ref{DicksonIdLemma}, the roots of the left side are
$w\ang{\zeta u}$ for $\zeta\in\mu_{q+1}$ and $w\in\F_q^\x$, and
these are distinct.  Denote the set of these roots by $S$; we have $|S|=q^2-1$,
and also $s\in S$ implies $ws\in S$ for all $w\in\F_q^\x$.
If $s_1,s_2\in S$, then $s_1^{q-1}=s_2^{q-1}$ iff 
$s_1/s_2\in\F_q^\x$. Thus, each power $s^{q-1}$ has exactly $q-1$
preimages in $S$, namely $\{ws : w \in\F_q^\x\}$. It follows that 
there are exactly $q+1$ distinct values $\{s^{q-1}:s\in S\}=\{\ang{\zeta u}^{q-1}:
\zeta\in\mu_{q+1}\}$.  This shows that the values $\ang{\zeta u}^{q-1}$ are distinct.
Since $R=e T(e)^2 X^{q-1}$, it follows that $e T(e)^2 \ang{\zeta u}^{q-1}$ 
are roots of $R^{q+1}+a^2 R + a^2$, and since they are distinct, all $q+1$ roots are accounted for.
\end{proof}

\section{ Splitting field of $x^{q+1} - b x + b$.} \label{qplus1Sec}
For this section only, we will consider both even and odd characteristic. Let $p$ be a
prime and $q=p^n$, where $n\ge 1$.
The polynomial $f(x) = x^{q+1}-bx+b$ in characteristic~$p$ (where $b\ne 0$) was studied in \cite{qplus1}.
The polynomial $x^{q+1}+a^2x+a^2$ that arises in Proposition~\ref{eRelation} is the
special case $p=2$, $b=a^2$.  The article \cite{qplus1} gives explicit formulas
for the splitting field and for the Galois action, which we recall in Theorems~\ref{qplus1Thm}
and~\ref{qplus1Thm2} below. Theorem~\ref{qplus1Thm} is illustrated in Figure~\ref{splitting_diagram}.

\begin{theorem} \label{qplus1Thm} (\cite{qplus1}) Let $q=p^n$ and let $b$ be a 
nonzero element of a field $\F$ in char.~$p$.
Let $r$, $r_0$, $r_1$ be distinct roots of
$x^{q+1}-bx+b$, and define
\begin{equation} \label{zyxi}
z=r_0/r,\qquad
y=(r_1-r)/(r_1-r_0),\qquad \xi=y^q-y.
\end{equation}
Then
$$y^{q-1}=z,\qquad \xi=y(z-1),\quad {\rm and}\quad
z(z-1)^{q-1}=1/(r-1)=\xi^{q-1}.$$
We have $y^{q^2}-y=\xi^q+\xi\ne 0$.
\end{theorem}

\begin{figure} 
\begin{equation*}
\xymatrix{
& \F(y)=\F(r,r_0,r_1) \ar@{-}[dl]^{q-1} \ar@{-}[dr]_q &\\
\F(z)=\F(r,r_0) \ar@{-}[dr]^{q}& &
\F(\xi) \ar@{-}[dl]_{q-1} \\
&\F(r) \ar@{-}[d]^{q+1} & \\
&\F&}
\end{equation*}
\caption{Splitting field of $x^{q+1}-bx+b$ when $\F=\F_q(b)$ and $b$ is transcendental. Here $r,r_0,r_1$ are roots of $f$,
$y=(r_1-r)/(r_1-r_0)$,
$z=y^{q-1}$, $\xi=y^q-y$, $\xi^{q-1}=z(z-1)^{q-1}=1/(r-1)$, $r_0=rz$,
$r_1=r(y-1)^{q-1}$. The degree of each
extension field is indicated. The minimal polynomial for $y$ over $\F(z)$ is $y^{q-1}-z$. The minimal polynomial
for $y$ over $\F(\xi)$ is $y^q-y-\xi$. The minimal polynomial for $z$ over $\F(r)$ is $z(z-1)^{q-1}-1/(r-1)$.
The minimal polynomial for $\xi$ over $\F(r)$ is $\xi^{q-1}-1/(r-1)$. The minimal polynomial for $r$ over~$\F$ is
$r^{q+1}-br+b$. \label{splitting_diagram}}
\end{figure}
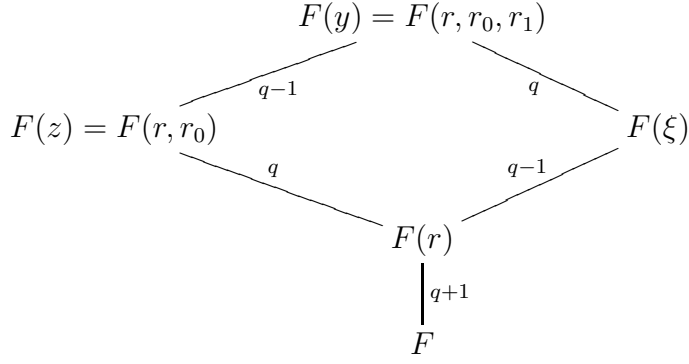

\begin{theorem} \label{qplus1Thm2}  Let $f(x)=x^{q+1}-bx+b$, and let $r$,
$r_0$, $r_1$, $y$ be as in Theorem~\ref{qplus1Thm}. Then $r_0=ry^q$ and
	$r_1=r(y-1)^{q-1}$.
The complete set of roots of $f$ is
$\{r_w : w \in \P^1(q)\}$, where $r_\infty=r$ and $r_w=r(y-w)^{q-1}$
for $w\in\F_{q}$. The roots are distinct. The splitting field over
$\F_p(b)$ is $\F_{q}(y)$. If $\sigma\in\Gal(f/\F_p(b))$,
then there is a unique $\gamma\in \PGL_2(q)$ such that 
$\sigma(y)=\gamma^{-1}(y)$. We have $\sigma(r_w)=r_{\gamma(\sigma w)}$,
where $\gamma$ has the usual action by linear fractional transformations
on $\P^1(q)$. For any $\gamma \in \PGL_2(\F_q)$,
\begin{equation}
	\gamma^{-1} (y) = \frac {r_{\gamma(1)} - r_{\gamma(\infty)}}{ r_{\gamma(1)} - r_{\gamma(0)}}. \label{gammay}
\end{equation}
\end{theorem}

\begin{proof}  All the above results were proved in \cite{qplus1} except for (\ref{gammay}), which we will
prove here.  Write $\g=\textmatrix abcd$,
so $\g^{-1}(y) = (dy-b)/(-cy+a)$.  Eq.~(\ref{gammay}) 
is equivalent to 
$$(dy-b)(r_{\g(1)}-r_{\g(0)}) = (-cy+a)(r_{\g(1)}-r_{\g(\infty)})$$
which in turn is equivalent to
\begin{equation} 
((c+d)y-(a+b))r_{\g(1)} = (dy-b)r_{\g(0)}+(cy-a)r_{\g(\infty)}. 
	\label{gammay2} 
\end{equation}

First assume that $\g(\infty)$, $\g(0)$, and $\g(1)$ are all finite.
Then $cd(c+d)\ne 0$.
Since $r_w=r(y-w)^{q-1}$ and $r\ne 0$, 
\begin{eqnarray*} \text{(\ref{gammay}) holds } &\iff& 
	\text{(\ref{gammay2}) holds} \\
&\iff& ((c+d)y-(a+b))(y-(a+b)/(c+d))^{q-1} = \\
	&&\quad (dy-b)(y-b/d)^{q-1} + (cy-a)(y-a/c)^{q-1} \\
&\iff& (c+d)(y-(a+b)/(c+d))^q = d(y-b/d)^q + c(y-a/c)^q\\
&\iff& (c+d)y^q-(a+b) = dy^q-b + cy^q-a.
\end{eqnarray*}
The final equation clearly is true, therefore (\ref{gammay}) holds
in this case.

The other three cases to consider are: $\g(\infty)=\infty$,
$\gamma(1)=\infty$, and $\g(0)=\infty$. 

If $\g(\infty)=\infty$, then $\gamma = \textmatrix ab01$, and 
eq.~(\ref{gammay2}) becomes 
$$(y-(a+b))r_{a+b}=(y-b)r_b -ar.$$
Since $r\ne0$ and $r_w/r=(y-w)^{q-1}$, this is equivalent to
$(y-(a+b))^q = (y-b)^q-a$, which is clearly true. Thus, (\ref{gammay})
holds when $\g(\infty)=\infty$.

If $\g(1)=\infty$, then $c+d=0$ and one can normalize to make
$c=-1$, $d=1$.  Eq.~(\ref{gammay2}) becomes 
	$$-(a+b) r = (y-b) r_{b} + (-y-a)r_{-a}$$
which is equivalent to $-(a+b)=(y-b)^q-(y+a)^q$.  The latter is true, so 
(\ref{gammay}) again holds.

Finally, if $\g(0)=\infty$, then $\gamma = \textmatrix ab10$ and
eq~(\ref{gammay2}) becomes
$$ (y-(a+b)) r_{a+b}= -br + (y-a)r_a.$$
This is equivalent to the true statement: $(y-(a+b))^q=-b+(y-a)^q$,
and the proof is complete.
\end{proof}

\noindent{\bf Remark.\ } If $w\in\F_q$, then $w$ can explicitly be expressed in terms of the
roots of $x^{q+1}-bx+b$ as follows. Let $\gamma = \textmatrix 1{ w} 01$. Then
\begin{eqnarray*} w &=& y - (y-w) = y - \gamma^{-1} (y) \\
	&=& \frac{r_1-r}{r_1-r_0} - \frac{r_{\gamma (1)}-r_{\gamma(\infty)}}{r_{\gamma (1)}-r_{\gamma (0)}} \\
&=&\frac{r_1-r}{r_1-r_0} - \frac{r_{w+1}-r}{r_{w+1}-r_w}.
\end{eqnarray*}

\begin{lemma} \label{bpLemma} Let $p$ be any prime (even or odd), $q=p^n$, 
$\F$ a field of characteristic~$p$ (not necessarily perfect),
and $0\ne b \in \F$. Then $x^{q+1}-bx+b$ and $x^{q+1} - x + 1/b$ have the same splitting field over~$\F$.
Also, the polynomials $x^{q+1} - b^{p^i} x + b^{p^i}$ have the same splitting field over $\F$ for all $i\ge 0$.
\end{lemma}

\begin{proof}  We begin by proving that $x^{q+1} -b x + b$ and $x^{q+1} - b^p x + b^p$ have the same splitting field
over~$\F$.  Denote these splitting fields by $\LL$ and $\LL_1$, respectively.  
Let $r_1,\ldots,r_{q+1}$ denote the roots of $x^{q+1}-bx+b$,
so that $\LL=\F(r_1,\ldots,r_{q+1})$ and $\LL_1=\F(r_1^p,\ldots,r_{q+1}^p)$.
	To show $\LL=\LL_1$, it suffices to demonstrate that $r_i\in \LL_1$.
	Indeed, $r_i^{q+1}=b(r_i-1)$ implies $1-1/r_i=r_i^q/b\in \LL_1$, therefore
	$r_i\in \LL_1$ as required.

We showed that $x^{q+1} - b^p x + b^p$ has the same splitting field over~$\F$ as $x^{q+1} - b x + b$.
Repeating the argument with $b^p$ in place of $b$, we see that $x^{q+1}-b^{p^2} x + b^{p^2}$
has the same splitting field over~$\F$ as $x^{q+1} - b^p x + b^p$. By induction on~$i$,
all fields $x^{q+1}-b^{p^i} x + b^{p^i}$ have the same splitting field over~$\F$.

It remains to prove that $x^{q+1} - x + 1/b$ has the same splitting field as well.
If $r$ is a root of $x^{q+1}-x+1/b$, then $br$ is a root of $x^{q+1}-b^q x + b^q$, because
$$(br)^{q+1} - b^q (br) + b^q = b^{q+1} (r^{q+1} - r + 1/b) = 0.$$
This shows $\SF(x^{q+1}-x+1/b;\F) = \SF(x^{q+1}-b^q x + b^q;\F)$. 
\end{proof}

\section{Expressing roots of $C(x)+a$ in terms of $\SF(x^{q+1}+ax+a;\F)$} \label{sec:RootsOfC}

Now we apply the theory from Section~\ref{qplus1Sec} to derive formulas expressing the roots of $C(x)+a$
in terms of the roots of $x^{q+1}+ax+a$, where $0\ne a \in \F$.  
For the remainder of this article, char$(\F)=2$ and $q=2^n>2$.   
(See {\it Notation} in Section~\ref{sec:Intro}.)
 
Let $e,u\in\cj \F$ satisfy
$$C(e)=a,\qquad 1/e = \ang{u^{q+1}}.$$
Proposition~\ref{eRelation} showed that the roots of $x^{q+1} + a^2 x + a^2$ are
$$\text{$\{\,\lambda \ang{\zeta u}^{q-1} : \zeta \in \mu_{q+1}\}$, \quad where
$\lambda = e\, T(e)^2$.}$$

Let $r$, $r_0$, $r_1$ be any three distinct roots of $x^{q+1} + a x + a$.  Then 
$r^2$, $r_0^2$, and $r_1^2$ are distinct roots of $x^{q+1} + a^2 x + a^2$. 
After rescaling $u$ by an element of $\mu_{q+1}$, we can arrange that $r^2 = \lambda \ang u^{q-1}$,
while still keeping the condition $1/e=\ang{u^{q+1}}$.
Then there are distinct $\zeta,\rho \in \mu_{q+1} \setminus \{1\}$ such that
\begin{equation*}
r^2 = \lambda \ang{u}^{q-1} \qquad
r_0^2 = \lambda \ang{\zeta^2 u}^{q-1} \qquad
r_1^2 = \lambda \ang{\rho^2 u }^{q-1}.
\end{equation*}
(Here we have used the fact that $\z\mapsto \z^2$ is a bijection of $\mu_{q+1}$
to itself.) Let 
$$y = (r_1-r)/(r_1-r_0).$$  
By Theorem~\ref{qplus1Thm},  the splitting field of $x^{q+1}+ax+a$ is 
$\LL=\F \circ \F_q(y)$.

\begin{lemma} \label{sevenFormulas} Let $y$, $e$, $\zeta$, $\rho$ be as above, and let
\begin{equation} c = \frac{\ang{\zeta/\rho}}{\ang\zeta\ang\rho}\qquad{\rm and}\qquad d = \frac 1 {\ang\zeta}. \label{cDef} \end{equation}
	Then $c\in\F_q^\x$, $d\in\F_{q,1}$, and the following formulas hold.
%\begin{equation} y^2 = \frac{\ang{\rho^2} \ang{\zeta^2 u} } {\ang{\zeta^2/\rho^2} \ang u } \label{formula1} \end{equation}
	\begin{equation} y = \frac { (d/c) (\zeta u + 1/\zeta) } {u + 1 } \label{formula2} \end{equation}
%\begin{equation} y = \frac { \ang{\rho} (\zeta u + 1/\zeta) } {\ang{\zeta/\rho} (u + 1) } \label{formula2} \end{equation}
		\begin{equation} u = \frac {(c/d) y + 1/\zeta } {(c/d) y + \zeta }  \label{formula3} \end{equation}
%\begin{equation} u = \frac {\ang{\zeta/\rho} y + \ang{\rho}/\zeta } {\ang{\zeta/\rho} y + \ang{\rho}\zeta }  \label{formula3} \end{equation}
\begin{equation} \ang{u} = \frac {1}{(cy)^2 + c y + d^2} \label{formula4} \end{equation}
\begin{equation} 1/e = D_{q+1}\left(\frac{1}{(cy)^2+cy+d^2}\right) \label{formula5} \end{equation}
\begin{equation}  (y^q+y)^2 = \frac { \ang {u^{q+1}} } {c^2{\ang u}^{q+1} } 
\label{formula6} \end{equation}
\begin{equation} e = \left(cy^2+y+\frac {d^2}{c} \right)^{q+1} \cdot (y^q+y)^{-2}. \label{formula7} \end{equation}
\begin{equation} c(y^q+y) = 1+\sum_{i=0}^{n-1} \ang u^{-2^i}. \label{formula8}
\end{equation}
\end{lemma}

\begin{proof}  For any $\l\in\mu_{q+1}\setminus \{1\}$, $\ang{\l}^q
	=\ang{\l^q}=\ang{\l^{-1}}=\ang{\l}$ and $\ang\lambda\ne 0$,
	so $\ang\lambda\in\F_q^\x$.
	Since 1, $\rho$, $\z$ are
	distinct, each of $\z,\rho,\z/\rho$ belongs to
	$\mu_{q+1}\setminus \{1\}$, therefore 
	$c,d\in\F_q^\x$. That $d\in\F_{q,1}$ was proved
	in Lemma~\ref{Fq1Lemma}.
	Now we compute:
\begin{eqnarray*} y^2 &=& \frac{r_1^2-r^2}{r_1^2-r_0^2}  \\
&=& \frac{\ang{\rho^2 u}^{q-1} - \ang{u}^{q-1}}{\ang{\rho^2 u}^{q-1} - \ang{\zeta^2 u}^{q-1}}  \times
\frac{ \ang{u} \ang{\rho^2 u} \ang{\zeta^2 u} }{ \ang{u} \ang{\rho^2 u} \ang{\zeta^2 u}} \\
&=& \frac{ \ang{\rho^2 u}^q \ang{u} + \ang{u}^q \ang{\rho^2 u} } {\ang{\rho^2 u}^q \ang{\zeta^2 u} + \ang{\zeta^2 u}^q \ang{\rho^2 u} }
\times \frac {\ang{\zeta^2 u}} {\ang{u}} \\
&=& \frac{ \ang{\rho^{-2} u^q} \ang{u} + \ang{u^q} \ang{\rho^2 u} } {\ang{\rho^{-2} u^q} \ang{\zeta^2 u} + \ang{\zeta^{-2} u^q} \ang{\rho^2 u} }
\times \frac {\ang{\zeta^2 u}} {\ang{u}}. \\
\end{eqnarray*}
In the first fraction, the numerator and denominator can be rewritten 
using~(\ref{angProperty}) as
\begin{eqnarray*} \ang{\rho^{-2} u^q} \ang{u} + \ang{u^q} \ang{\rho^2 u}  
%&=& (\rho^{-2} u^q + \rho^2 u^{-q}) (u+1/u) + (u^q + u^{-q})(\rho^2 u + \rho^{-2}u^{-1}) \\
&=& \ang{\rho^{2}} \ang{ u^{q+1}};  \\
	\ang{\rho^{-2} u^q} \ang{\zeta^2 u} + \ang{\zeta^{-2} u^q} \ang{\rho^2 u} 
&=& \ang{\zeta^2/\rho^2} \ang{u^{q+1}} = (c\ang\zeta\ang\rho)^2 \ang{u^{q+1}}
	= (c/d)^2 \ang{\rho^2} \ang{u^{q+1}}.
\end{eqnarray*}
Hence,
\begin{equation} y^2 = \frac{d^2\ang{\zeta^2 u} } {c^2 \ang u }. \label{formula1} \end{equation}
Multiply the right side of (\ref{formula1}) by 
	$u/u$ to obtain $y^2=(d^2/c^2)(\z^2u^2+\z^{-2})/(u^2+1)$.  Take the square root to obtain (\ref{formula2}).
Now (\ref{formula2}) shows that $u$ and $y$ are related by a linear fractional transformation over $\F_{q^2}$. Solving
for $u$ in terms of $y$ gives~(\ref{formula3}). To derive (\ref{formula4}), let 
$A=c y + d/\zeta$,
$B=cy+d\zeta$, so $u=A/B$. Then $\ang u = A/B+B/A=(A^2+B^2)/(AB)$.
Since $A^2+B^2 = d^2\ang{\zeta^2}=1$ and 
	$AB=c^2y^2+\ang\zeta cdy + d^2= c^2y^2+cy+d^2$, formula~(\ref{formula4}) follows.
We have $1/e=\ang{u^{q+1}}=D_{q+1}(\ang{u})$, which proves~(\ref{formula5}).
For~(\ref{formula6}), 
	$$ y = (d/c) \frac{\zeta u + 1/\zeta}{u+1},\qquad y^q = (d/c) \frac{\zeta^{-1}u^q+\zeta}{u^q+1},$$
therefore
\begin{eqnarray*} 
	y^q+y&=&(d/c) \cdot \frac{ (\zeta u + 1/\zeta)(u^q+1) + (\zeta^{-1} u^q + \zeta)(u+1)} {(u^q+1)(u+1) } \\
	&=& \frac{(d/c) \ang{\zeta}(u^{q+1}+1)}{(u+1)^{q+1}} = 
\frac{(u^{q+1}+1)}{c(u+1)^{q+1}}. 
\end{eqnarray*}
Now square both sides and multiply on the right by $u^{-(q+1)}/u^{-(q+1)}$ to obtain~(\ref{formula6}). 
Formula (\ref{formula7}) is obtained by substituting $1/e=\ang{u^{q+1}}$ and $1/\ang u = c^2 y^2 + c y + d^2$
into~(\ref{formula6}). By (\ref{formula6}) and Lemma~\ref{lem:Dicksonqm1},
	\begin{equation*} c^2(y^q+y)^2 \ang u^{q+1} = \ang{u^{q+1}}
		= D_{q+1}(\ang u) = \ang u^{q+1} + \sum_{i=1}^n \ang u^{q-2^i+1}.
\end{equation*}
Divide by $\ang u^{q+1}$ then take the square root to obtain~(\ref{formula8}).
\end{proof}

\begin{figure} [t]
\begin{equation*}
\xymatrix{
&& \F(y)=\F(r,r_0,r_1) \ar@{-}[ddll]^{q-1} \ar@{-}[d]_2 \ar@{-}[ddr]^{2(q+1)}&\\
&& \F(\ang{u}) \ar@{-}[d]_{q/2} \ar@{-}[dr]_{q+1} & \\
\F(z)=\F(r,r_0) \ar@{-}[dr]^{q} & &
\F(\xi) \ar@{-}[dl]_{q-1} &
\F(e) \ar@{-}[ddl]^{(q/2)(q-1)} \\
&\F(r) \ar@{-}[dr]^{q+1} && \\
&&\F=\F_q(a)&}
\end{equation*}
\caption{Joint splitting field of $x^{q+1}+ax+a$ and $C(x)+a$ when $q=2^n>2$, $\F=\F_q(a)$ and $a$ is transcendental. Here 
	$r,r_0,r_1$ are any distinct roots of 
$x^{q+1}+ax+a$ and $e$ is an arbitrary root of $C(x)+a$. As before,
$y=(r_1-r)/(r_1-r_0)$,
$z=r_0/r=y^{q-1}$, $\xi=y^q-y$, and $\xi^{q-1}=z(z-1)^{q-1}=1/(r-1)$.  
The formulas from Lemma~\ref{sevenFormulas} show that the minimal polynomial for $y$ over $\F(\ang{u})$ is 
$(cy)^2+cy+d^2-\ang{u}^{-1}$, where $c\in\F_q^\x$ and $d\in \F_{q,1}$ depend upon the choice of the root of $C(x)$.
The minimal polynomial for $\ang{u}$ over $\F(e)$ is $D_{q+1}(\ang u)-1/e$.  The minimal polynomial for $y$ over $\F(e)$ is 
$(cy^2+y+d^2/c)^{q+1}+e(y^q+y)^2$. 
The minimal polynomial for $1/\ang u$ over $\F(\xi)$ is
$\sum_{i=0}^{n-1} (1/\ang u)^{2^i} + c\xi + 1$.\label{splitting_diagram2}}
\end{figure}
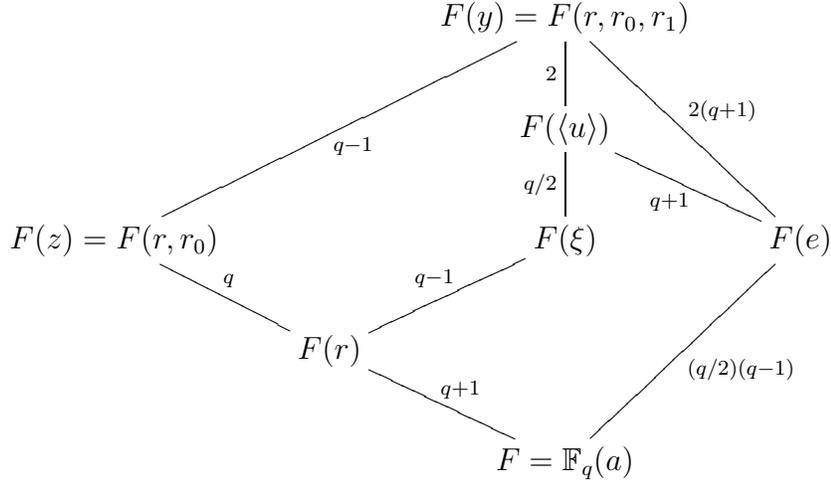

On account of Lemma~\ref{sevenFormulas}, Figure~\ref{splitting_diagram} can be extended to incorporate
other subfields of the splitting field, as shown in Figure~\ref{splitting_diagram2}.

\begin{theorem} \label{KsubsetL} Let $y=(r_1-r)/(r_1-r_0)$, where $r,r_0,r_1$ 
are three distinct roots of $x^{q+1}+ax+a$.
Then the complete set of roots of $C(x)+a$ is
\begin{equation*}
{\cal E} = \left\{\, \frac{ ( cy^2 + y + j/c)^{q+1} }{(y^q + y )^2} 
	: c \in \F_q^\x,\ j \in \F_{q,1} \, \right\}.
\end{equation*}
These roots are distinct.
If $X$ is an indeterminate, then
\begin{equation}\prod_{c \in \F_q^\x,\ j \in \F_{q,1}} \left(X - \frac{ (c y^2 + y + j/c)^{q+1} }{(y^q+y)^2}\right) = C(X) + a.\label{CXaIdentity}\end{equation}
\end{theorem}

\begin{proof} By Lemma~\ref{sevenFormulas} (in particular, (\ref{formula7})), 
if $C(e)=a$, then $e \in {\cal E}$. Thus, the roots of $C(x)+a$
are contained in ${\cal E}$.
The number of pairs $(c,j)$ is $(q-1)q/2$, which is the same
as $\deg(C)$. The roots of
$C(x)+a$ are distinct by Lemma~\ref{lem:distinctRootsC}. By the pigeon-hole
principle, ${\cal E}$ is the full set of roots and
the pairs $(c,j)$ give rise to distinct roots.
The last sentence follows from 
$C(X)+a = \prod_{e \in {\cal E}} (X-e)$. 
\end{proof}

\begin{corollary}  \label{cor:yIdentity} 
The following identity holds for all $y$: 
\begin{equation}
\prod_{c\in\F_q^\x,\ j \in \F_{q,1} } (c y^2 + y + j/c) = 1 + (y^q+y)^{q-1}.
	\label{borelIdentity}
\end{equation}
\end{corollary}

\begin{proof} We express $a$ in two ways.
First, substituting $X=0$ into (\ref{CXaIdentity}) gives
\begin{equation}\prod_{c \in \F_q^\x,\ j \in \F_{q,1}} \frac{ (c y^2 + y + j/c)^{q+1} }{(y^q+y)^2}= a.\label{aIdentity}\end{equation}
Second, by Theorem~\ref{qplus1Thm}, $\xi = y^q-y$,
$\xi^{q-1}=1/(r+1)$, and $r^{q+1}+a r + a = 0$. Hence,
$$a=\frac{r^{q+1}}{r+1} = r^{q+1}\xi^{q-1} = (1+\xi^{1-q})^{q+1}\xi^{q-1},
\quad{\rm where}\quad \xi = y^q-y.$$
Comparing the two expressions, we find that
	$$\frac{\prod_{c\in\F_q^\x,\ j\in\F_{q,1}} (c y^2 + y + j/c)^{q+1}} { \xi^{q(q-1)} }  = (1+\xi^{1-q})^{q+1} \xi^{q-1}.$$
If $a$ (and hence also $y$) is transcendental, then this may be interpreted as an identity in the ring $\F_q(y)$.
On multiplying through by $\xi^{q(q-1)}$ and then taking the unique $(q+1)$th root belonging to $\F_q(y)$, we obtain:
$$\prod_{c\in\F_q^\x,\ j \in \F_{q,1}} (cy^2 + y + j/c) = \xi^{q-1}(1+\xi^{1-q}) =1+\xi^{q-1} = 1+(y^q+y)^{q-1}.$$
\noindent
\end{proof}

\noindent{\it Note added in proof:} Corollary~\ref{cor:yIdentity} can
also be proved using \cite[Proposition~3.4]{ArtinMaps}, as follows. Let
$L$ and $R$ denote the left and right sides of~(\ref{borelIdentity}).
Both $L$ and $R$ are easily seen to be invariant 
under the Borel group $\{\textmatrix ab01 : a\in\F_q^\x,\ b\in\F_q\}$, 
and their degrees are equal to the order of 
the Borel group.  By \cite[Proposition~3.4]{ArtinMaps}, 
$L=aR+b$ where $a\in\F_q^\x$ and $b\in\F_q$. 
Since $L$ and $R$ are monic, $a=1$, 
and since $L$ and $R$ have the same constant term, $b=0$.

%\begin{proposition} For $y\in\calY$, define $e(y,c,j)=(cy^2+y+j/c)^{q+1}/(y^q+y)^2$, so that $\{e(y,c,j) : c\in\F_q^\x, j \in \F_{q,1}\}$ is a complete
	%set of roots of $C(x)+a$. Let $d=j^{1/2}$, and let $\zeta\in\mu_{q+1}\setminus \{1\}$ such that $1/d=\ang\zeta$. Let
	%$$M = \begin{pmatrix} c/d & 1/\zeta \\ c/d & \zeta \end{pmatrix}.$$
		%Then
		%$$1/e(y,j,c)=\ang{u^{q+1},\quad{\it where}\ u=M(y).$$
%\end{proposition}
%\begin{proof} Let $e=e(y,c,j)$.

\section{Equality of splitting fields} \label{sec:Equality}
In this section, we prove one of our main results, that $x^{q+1}+ax+a$ and $C(x)+a$ have the same splitting field. This will be accomplished
by explicitly writing the roots of each polynomial in terms of the roots of the other.

From here on, let
\begin{equation} {\cal Y}
=\left\{\frac{r_1-r}{r_1-r_0} : r,r_0,r_1\text{\ are distinct roots of 
$x^{q+1}+ax+a$} \right\}. \label{calY}\end{equation}
For $y\in {\cal Y}$, $c\in \F_q^\x$ and $j\in \F_{q,1}$, define
\begin{equation} e(y,c,j) = 
\frac{ ( cy^2 + y + j/c)^{q+1} }{(y^q - y )^2}. \label{eycdDef} \end{equation}
By Theorem~\ref{KsubsetL}, for a fixed $y$, the values $\{ e(y,c,j) : c \in \F_q^\x,\ j \in \F_{q,1}\}$ are distinct and comprise the complete set of 
roots of $C(x)+a$. Note that $cY^2+Y + j/c$ comprise all irreducible
quadratics in $\F_q[Y]$, normalized with the $Y$ coefficient equal to unity.

\begin{theorem} \label{KeqL} For $q=2^n \ge 4$, $\SF(C(x) + a; \F) = \SF(x^{q+1} + x + 1/a;\F)= \SF(x^{q+1}+ax+a;\F)$.
\end{theorem}

\begin{proof} The equality $\SF(x^{q+1}+x+1/a;\F)=\SF(x^{q+1}+ax+a;\F)$ was shown in Proposition~\ref{bpLemma}, so it suffices
to show that $\KK=\LL$, where
$$\KK=\SF(C(x)+a;\F),\qquad \LL = \SF(x^{q+1}+ax+a;\F).$$
Theorem~\ref{KsubsetL} explicitly expresses each root of $C(x)+a$ in terms of the roots of $x^{q+1} + ax + a$.
(See the remark following Theorem~\ref{qplus1Thm2} to see how $c,j$ can be written in terms
of roots of $x^{q+1}+ax+a$.)  This implies that $\KK\subset \LL$.  
To show $\LL\subset \KK$, we will express an arbitrary root of $x^{q+1}+ax+a$ in terms
of the roots of $C(x)+a$.

Let $r$ be an arbitrary root of $x^{q+1}+ax+a$.  Select any other two roots 
$r_0$ and $r_1$ and define $y=(r_1-r)/(r_1-r_0)\in {\cal Y}$, 
$\xi=y^q-y$.
By Theorem~\ref{qplus1Thm}, $r = 1+\xi^{1-q}$, so it suffices to express $\xi$ in terms of the roots of $C(x)+a$.  These
roots are $\{\,e(y,c,j) : c \in \F_q,\ j \in \F_{q,1}\,\}$.  We have
\begin{eqnarray*}
e(y,c,j) &=& (cy^{2q}+y^q+j/c)(cy^2+y+j/c)/\xi^2 \\ 
&=&  \frac{c^2 y^{2q+2} + c(y^{2q+1} + y^{q+2}) + y^{q+1} +j \xi^2 +  (j/c) \xi + (j/c)^2}{\xi^2}.
\end{eqnarray*}

First assume that $q>4$.  
We claim there are $c_1,c_2,c_3,c_4\in \F_q^\x$ such that $\sum_{i=1}^4 c_i = 1$ and $\sum_{i=1}^4 1/c_i = 0$.
To see this, select $\alpha \in \F_q \setminus \F_4$; such $\alpha$ exists because $q>4$. 
Let $w_1 = 1$, $w_2 = \alpha$, $w_3 = 1/\alpha$, and note that $w_1 + w_2 + w_3 = 1 + \alpha + \alpha^{-1} = 1/w_1 + 1/w_2 + 1/w_3$.  
Furthermore, this value does not belong to $\F_2$, because 
all solutions to $x + 1/x + 1 \in \{0,1\}$ belong to $\F_4$.  
Set $w_4=1/(1+\alpha+\alpha^{-1})$. Then $1/w_1 + 1/w_2 + 1/w_3 + 1/w_4 = 0$. Since $w_4 \not \in \F_2$, we know
$w_4 \ne w_4^{-1}$, and so $w_1 + w_2 + w_3 + w_4 = 1 + \alpha + \alpha^{-1} + w_4 = 1/w_4 + w_4 \ne 0$.  Setting $c_i = w_i/\sum w_i$, we find
that $\sum c_i = 1$ and $\sum c_i^{-1} = 0$. This establishes the claim.  Then
	$$ \sum_{i=1}^4 e(y,c/c_i,j) = \frac{(j/c)\xi + (j/c)^2}{\xi^2},$$
therefore 
\begin{equation} \frac b \xi + \frac {b^2}{\xi^2} \in \KK\quad
	\text{for all $b\in \F_q^\x$}. \label{bxi}
\end{equation}
	Next, $e(y,c,j)+\sum_{i=1}^4 e(y,cc_i,j) = (y^{q+1}+j\xi^2)/\xi^2 + (j/c)/\xi + (j/c)^2/\xi^2 \in \KK$.  By (\ref{bxi}), it follows that
$j+y^{q+1}/\xi^2\in \KK$ for all $j\in\F_{q,1}$, and consequently
	$$ j_1 + j_2 \in \KK\qquad \text{for any $j_1,j_2\in\F_{q,1}$.}$$

Select distinct values 
$j_1,j_2,j_3 \in \F_{q,1}$ such that $j_i + j_k \ne 1$ for
each $i,k$.  To see that these exist, note that if $n$ is odd, then $1\in \F_{q,1}$, and so the sum of two elements of $\F_{q,1}$ is never one
and it suffices to select $j_1,j_2,j_3$ to be distinct. Since $|\F_{q,1}| = q/2 \ge 4$, this selection is possible. If $n$ is even, then
$q/2 \ge 8$, so there are at least eight choices for $j_1\in \F_{q,1}$, six choices for $j_2\in\F_{q,1}\setminus \{j_1,j_1+1\}$,
and four choices for $j_3 \in \F_{q,1} \setminus \{j_1,j_2,j_1+1,j_2+1\}$.  This shows again that $j_1,j_2,j_3$ can be selected
so that $j_i+j_k \ne 1$ for each pair $(i,k)$. Let $\tau_1 = j_1 + j_2$ and $\tau_2 = j_1 + j_3$, and note that
$\{\tau_1,\tau_2,\tau_1+\tau_2\} \cap \F_2 = \emptyset$.  Let
$$w_1 = \frac 1 {\tau_1(\tau_1+\tau_2)},\qquad
w_2 = \frac 1 {\tau_2(\tau_1+\tau_2)}$$
and observe that
$$\tau_1 w_1 + \tau_2 w_2 = 0,\qquad \tau_1^2 w_1 + \tau_2^2 w_2 = 1.$$
Since $\tau_i$ and $w_i/\xi + w_i^2/\xi^2$ are in $\KK$, so is
$$\tau_1^2 (w_1/\xi + w_1^2/\xi^2) +\tau_2^2 (w_2/\xi + w_2^2/\xi^2) = 1/\xi.$$
This shows $\xi \in \KK$. Since $\xi^{q-1}=1/(1-r)$, $r\in \KK$. 
Since $r$ is an arbitrary root of $x^{q+1} + x + 1/a$, this completes the proof when $q>4$.

If $q=4$, then let $\alpha\in \F_4\setminus \F_2$. By direct calculation, 
$$1/\xi = e(y,1,\alpha) + e(y,1,\alpha^2) + \frac{ \left(e(y,\alpha,\alpha) + e(y,\alpha,\alpha^2) \right) 
\left(e(y,\alpha^2,\alpha) + e(y,\alpha^2,\alpha^2)\right) }
{e(y,1,\alpha) + e(y,1,\alpha^2)}.$$
	Thus, $\xi \in \KK$, and consequently $r\in \KK$, as desired.
\end{proof}

\section{Galois action and related factorizations} \label{galoisSec}

Since the roots of the two polynomials $C(x)+a$ and $x^{q+1} + a x + a$ 
generate the same field $\LL$,
any element of the Galois group $\Gal(\LL/\F)$ simultaneously permutes the roots of $C(x)+a$ and of $x^{q+1} + a x + a$.
For this reason, the factorizations of these two polynomials are related.  This section explores this. 

Let ${\cal Y}$ as in (\ref{calY}) and let $y\in{\cal Y}$, so 
$\LL=\F_q(y)$ is the splitting field of the two polynomials. 
Recall that the distinct roots of $C(x)+a$ are 
$\{\,e(y,c,j) : c \in \F_q^\x,\ j\in \F_{q,1} \}$, where
$e(y,c,j) = (c y^2 + y + j/c)^{q+1}/(y^q+y)^2$. 
If $\gamma \in \PGL_2(\F_q)$, then 
$\gamma^{-1}(y)\in {\cal Y}$ by (\ref{gammay}), therefore
$\{e(\gamma^{-1}(y),c,j): c\in\F_q^\x,\ j\in \F_{q,1}\}$ is also a complete
set of roots of $C(x)+a$. The next lemma shows how these two sets of roots
are related.

\begin{lemma} \label{galCLemma} If $b \in \F_q^\x$, then
\begin{eqnarray}
e(y+b,c,j) &=& e(y,c,j+bc+(bc)^2)  \label{yplusb} \\
e(by,c,j) &=& e(y,bc,j)  \label{yb} \\
e(1/y,c,j) &=& e(y,j/c,j). \label{recipy}
\end{eqnarray}
\end{lemma}

\begin{proof} For the first formula, 
\begin{eqnarray*} e(y+b,c,j) &=&
\frac{(c(y+b)^2+(y+b)+j/c)^{q+1}}{((y+b)^q+(y+b))^2} \\
&=& \frac{(cy^2+y + (j+bc+b^2c^2)/c)^{q+1}}{(y^q+y)^2} \\
&=& e(y,c,j+bc+(bc)^2).
\end{eqnarray*}
Note that $j+cb+c^2b^2 \in \F_{q,1}$.
For the second formula, 
$$e(by,c,j)=\frac{(cb^2y^2 + by + j/c)^{q+1}}{((by)^q+by)^2} = \frac{ b^{q+1} (cby^2+y+j/(cb))^{q+1} }
{ b^2 (y^q+y)^2 } = e(y,bc,j).$$
Finally,
\begin{eqnarray*}
e(1/y,c,j) &=& \frac{ (c y^{-2} + y^{-1} + j/c)^{q+1} } {(y^{-q}+y^{-1})^2 } \times \frac{ y^{2(q+1)} }{ y^{2(q+1)} }
= \frac{(c+y+(j/c) y^2)^{q+1}} {(y+y^q)^2} = e(y,j/c,j).
\end{eqnarray*}
\end{proof}

\begin{theorem} \label{sigmaFacThm} Let $\LL$ be the splitting field of $x^{q+1}+x+1/a$ (which is also the splitting
field of $C(x)+a$ by Theorem~\ref{KeqL}), and let $\sigma \in \Gal(\LL/\F\circ \F_q)$.
\begin{itemize}
\item[({\it i})] If $\sigma$ fixes at least three roots of $x^{q+1} + x + 1/a$, then it fixes all roots of $x^{q+1}+x+1/a$
and all roots of $C(x)+a$.
\item[({\it ii})] If $\sigma$ fixes exactly two roots of $x^{q+1}+x+1/a$, 
then the permutation induced by $\sigma$ on the roots of $x^{q+1}+x+1/a$ 
has orbits of size $1,1,\delta,\delta,\ldots,\delta$, where $\delta$ 
divides $q-1$. The permutation induced by $\sigma$ 
on the roots of $C(x)+a$ has all its orbits of size $\delta$.
\item[({\it iii})] If $\sigma$ fixes exactly one root of $x^{q+1} + x + 1/a$, then the remaining roots fall into
$\sigma$-orbits of size~2.
Also, $\sigma$ fixes exactly $q/2$ roots of $C(x)+a$, and the remaining roots fall into exactly $q/2(q/2-1)$ $\sigma$-orbits
of size two. 
\item[({\it iv})] If $\sigma$ fixes no roots of $x^{q+1} + x + 1/a$, then all $\sigma$-orbits of 
$x^{q+1}+x+1/a$ have the same size $\delta$, where $\delta$ divides $q+1$.
Also, $\sigma$ fixes exactly one root of $C(x)+a$, and all remaining roots belong to $\sigma$-orbits of size $\delta$.
\end{itemize}
\end{theorem}

\begin{proof} The two polynomials $x^{q+1}+x+1/a$ and $x^{q+1} + a x + a$ have the same splitting field by Lemma~\ref{bpLemma}, and their roots are in bijection
by the Galois-invariant map $r\mapsto ar^{1/q}$, therefore the cycle structure
of the permutation induced by $\sigma$ on the roots is the same for both 
polynomials.
For this reason, we can work with $x^{q+1}+ax+a$ instead of $x^{q+1}+x+1/a$.  

By Theorem~\ref{qplus1Thm2}, if $r,r_0,r_1$ are any three distinct roots of 
$x^{q+1}+ax+a$ and $y=(r_1-r)/(r_1-r_0)$, then $\LL = \F \circ \F_q(y)$, and
there is a unique $\g\in \PGL_2(\F_q)$ such that $\sigma(y) = \gamma^{-1}(y)$.
Furthermore, $\sigma(r_w)= r_{\gamma(\sigma(w))}$, where $r_w$ for 
$w\in\P^1(q)$ are 
defined in Theorem~\ref{qplus1Thm2} and comprise a complete set of roots 
of $x^{q+1}+ax+a$. (In this labeling, $r=r_\infty$.)
Since by hypothesis $\s$ fixes $\F_q$, $\sigma(r_w)=r_{\gamma(w)}$.

In case ({\it i}), let $r,r_0,r_1$ be three distinct roots of $x^{q+1}+ax+a$ that are fixed by $\sigma$. Then 
$y=(r_1-r)/(r_1-r_0) $ is also fixed by $\sigma$.  By Theorem~\ref{qplus1Thm}, $\LL=\F\circ \F_q(y)$.
Since $\sigma\in\Gal(\LL/\F\circ \F_q)$ and it fixes $y$, it follows that $\sigma$ is the identity, and so it fixes
all roots of both polynomials.

In case ({\it ii}), let $r$ and $r_0$ be roots of $x^{q+1}+ax+a$ that are fixed by $\sigma$, and 
select a third root $r_1$ with which to form $y$.
	Since $r=r_\infty$ and $\sigma(r_w)=r_{\gamma(w)}$, 
$\gamma$ fixes $\infty$ and~0. Then it has the form $\textmatrix b001$ with $b\in\F_q^\x$.
Let $\delta$ be the multiplicative order of $b$. The orbits of $\gamma$ acting on $\P^1(q) \setminus \{\infty,0\}$
are of the form $\{w,bw,b^2w,\ldots, b^{\delta-1}w\}$, showing that the non-singleton orbits all have the same order $\delta$.
Consequently, the $\sigma$-orbits on roots of $x^{q+1}+x+1/a$ have sizes $1,1,\delta,\delta,\ldots,\delta$.
The action on the roots of $C(x)+a$ is
	$$\sigma\left( e(y,c,j)\right)=e(\gamma^{-1}(y),c,j)=
	e(b^{-1}y,c,j) = e(y,c/b,j) \qquad (c\in\F_q^\x,\ j \in \F_{q,1}).$$
Thus, each $\sigma$-orbit has size exactly $\delta$.

In case~({\it iii}), $\sigma$ fixes exactly one root of $x^{q+1}+ax+a$, which we call $r$. 
Let $r_0,r_1$ be two other roots and define $y=(r_1-r)/(r_1-r_0)$. 
Then $\sigma(y)=\gamma^{-1}(y)$, where $\gamma(\infty)=\infty$ and
$\gamma$ fixes no other element of $\P^1(\F_q)$.
The elements of $\PGL_2(\F_q)$ that fix only $\infty$ are
$\g=\textmatrix 1b01$ with $b\in\F_q^\x$. 
Then order$(\sigma)={\rm order}(\gamma)=2$, $\sigma(y)=y+b$, and
$$\sigma\left( e(y,c,j)\right)=e(y+b,c,j)=e(y,c,j+(bc)+(bc)^2).$$
Since $b$ and $c$ are nonzero, $(c,j)=(c,j+(bc)+(bc)^2)$ iff 
$c=1/b$. So the $q/2$ roots $e(y,1/b,j)$ with $j\in\F_{q,1}$ are fixed, 
and the other roots belong to $\sigma$-orbits of size~2.

In case~({\it iv}), define $y$ with respect to three roots $r$, $r_0$, $r_1$ of $x^{q+1}+ax+a$ such that $\sigma(r) = r_0$ and
$\sigma(r_0) = r_1$.  
Since $\gamma$ takes $\infty$ to 0
and 0 to 1, it has the form $\gamma = \textmatrix 0 1 k 1$.  By hypothesis, 
$\sigma$ fixes no roots, and so
$1/(kw+1)=w$ has no solutions in $\F_q$.  This is equivalent to $(kw)^2 + (kw) + k$ having no rational roots, which is equivalent
to $k \in \F_{q,1}$. Let $\delta={\rm order}(\gamma)$.
We will show in Proposition~\ref{prop:dihedral} (or see \cite{Dickson}) that 
$\delta$ divides $q+1$ and that $\gamma^i$ has
no fixed points in $\P^1(\F_q)$ for $1\le i < \delta$. (In the notation
of (\ref{DjCDef}), $\gamma$ belongs to the group ${\cal C}_{\sqrt k,\sqrt k}$, 
which is cyclic of order~$q+1$.)
Thus, the orbits of $\sigma$ on the roots of $x^{q+1}+ax+a$ all have the same size, $\delta$.

	We have $\gamma^{-1} (y) = (y+1)/(ky)$, and so 
	\begin{eqnarray*} \sigma(e(y,c,j)) &=& e((y+1)/(ky),c,j) \\
&=& e(1+1/y,c/k,j) \qquad \text{by (\ref{yb})} \\
&=& e(1/y,c/k,j+c/k+c^2/k^2)\qquad \text{by (\ref{yplusb})} \\
&=& e(y,jk/c+1+c/k,j+c/k+c^2/k^2)\qquad \text{by (\ref{recipy}).} 
\end{eqnarray*}
This can equal $e(y,c,j)$ only if $(c/k)+(c/k)^2=0$, {\it i.e.} $c=k$.  In that case, 
$$ \sigma(e(y,k,j))=e(y,j,j).$$
So, for the root to be fixed by $\sigma$, we also need $j=k$.  Thus, there is 
exactly one fixed root, namely $e(y,k,k)$.  The other roots must belong to orbits of size dividing $\delta$,
where $\delta$ is the order of $\gamma$. We claim the orbits have size exactly $\delta$.  Indeed, suppose that $\sigma$
had an orbit of size~$i$, where $i$ strictly divides $\delta$, and consider $\sigma^i$.
This fixes no roots of $x^{q+1}+ax+a$, so it fixes exactly one root of $C(x)+a$, which must be $e(y,k,k)$.  But it
also fixes the points on the $\sigma$-orbit of size~$i$, a contradiction.
So the roots of $C(x)+a$ fall into $\sigma$-orbits of size $1,\delta,\delta,\ldots,\delta$.
\end{proof}

\begin{corollary}  \label{facTypes}
For a polynomial $g\in \F[x]$, write $g \sim [n_1,n_2,...,n_t]$ if $g$ factors into irreducibles 
of degrees $n_1, n_2, \ldots , n_t$.  
Let $0\ne a\in \F=\F_{q^m}$, where $q=2^n>2$ and $m\ge1$. 
Then one of the following occurs.
\begin{itemize}
\item[({\it i})]  $x^{q+1}+x+1/a\sim [1,1,\ldots,1]$ and
	$C(x)+a\sim [1,1,\ldots,1]$;
\item[({\it ii})] 
	$x^{q+1}+x+1/a\sim[1,1,\delta,\delta,\ldots,\delta]$ and $C(x)+a\sim[\delta,\delta,\ldots,\delta]$, where $\delta|q-1$.
\item[({\it iii})] $x^{q+1}+x+1/a\sim[1,2,2,\ldots,2]$ and 
	$C(x)+a\sim[1,1\ldots,1,2,2,\ldots,2]$,
where $C(x)+a$ has $q/2$ linear factors and $(q^2-2q)/4$  irreducible quadratic factors.
\item[({\it iv})] 
	$x^{q+1}+x+1/a\sim[\delta,\delta,\ldots,\delta]$ and $C(x)+a\sim[1,\delta,\delta,\ldots,\delta]$,
where $\delta|q+1$.
\end{itemize}
\end{corollary}

\begin{proof}  Apply Theorem~\ref{sigmaFacThm}, taking $\sigma$ to be the Frobenius map: $\sigma(u) = u^{|\F|}$.  The sizes of the $\sigma$-orbits
	acting on the roots of $x^{q+1}+ax+1/a$ or $C(x)+a$ are the degrees of the irreducible factors over~$\F$.
\end{proof}

\begin{corollary} If $\F = \F_{q^m}$ and $0\ne a \in \F$, then $C(x)+a$ has exactly 0, 1, $q/2$, or $(q/2)(q-1)$ roots in~$\F$.
Let $c_i$ denote the number of $a \in \F^\x$ for which $C(x) + a$ has exactly $i$ roots in~$\F$. If $m$ is even, then
$$c_0 = \frac{(q-2)(q^m-1)}{2(q-1)},\qquad c_1 = \frac{q^{m+1}-q}{2(q+1)},\qquad c_{q/2} = q^{m-1},\qquad c_{(q/2)(q-1)} = \frac{q^{m-1}-q}{q^2-1}.$$
If $m$ is odd, then
$$c_0 = \frac{(q-2)(q^m-1)}{2(q-1)},\qquad c_1 = \frac{q^{m+1}+q}{2(q+1)},\qquad c_{q/2} = q^{m-1}-1,\qquad c_{(q/2)(q-1)} = \frac{q^{m-1}-1}{q^2-1}.$$
\end{corollary}

\begin{proof}  For $i\in \{0,1,2,q+1\}$, let $N_i$ denote the number of $a \in \F^\x$ such that $x^{q+1}+ax+a$ has exactly $i$ roots. By
Corollary~\ref{facTypes}, $N_0=c_1$, $N_1 = c_{q/2}$, $N_2 = c_{0}$, and $N_{q+1} = c_{(q/2)(q-1)}$.  The $N_i$'s are computed in
\cite[Theorem~5.6]{qplus1}. The result follows.
\end{proof}

We conclude this section by proving Proposition~\ref{quinticSextic}, which gives the related factorizations of $x^{q+1}+x+1/a$ and $C(x)+a$ when 
$q=4$ and $\F=\F_{2^k}$. We are not assuming
that $\F_q \subset \F$. The polynomials are
$x^5+x+1/a$ and $x(x+1)^5+a$.  Since $x^5+ax+a$ has the same splitting field and factorization type as $x^5+x+1/a$, we may study it instead.
Let $\LL$ denote the splitting field and let $\sigma \in \Gal(\LL/\F)$ denote the Frobenius element, $\sigma(b)=b^{|\F|}$; then $\sigma$ generates
$\Gal(\LL/\F)$.

If $k$ is even, then $\F_q\subset \F$. In that case, Proposition~\ref{quinticSextic} follows from Corollary~\ref{facTypes}. 

Now assume $k$ is odd,  and we must show that one of the following cases holds.
$$\text{$x^5+x+1/a \sim [1,1,1,2]$ and $x(x+1)^5 + a \sim [2,2,2]$}$$
$$\text{$x^5+x+1/a \sim [1,4]$ and $x(x+1)^5  + a\sim [1,1,4]$}$$
$$\text{$x^5+x+1/a \sim [2,3]$ and $x(x+1)^5  + a\sim [6]$}.$$
Note that $\sigma(c)=c^2$ for $c\in\F_4$.

If $\sigma$ fixes at least three roots of $x^{q+1}+ax+a$, then we can arrange for $y$ to be rational. 
Then $\sigma(r_w) = r_{\sigma(w)}$ for $w\in \P^1(\F_4)$. Let
$\a\in\F_4\setminus\F_2$. The conjugate pair
$\alpha$ and $\alpha^2$ are exchanged by $\sigma$,
while all other elements of $\P^1(q)$ are fixed, and so 
$x^5+ax+a \sim [1,1,1,2]$.
The roots of $C(x)+a$ are $e(y,c,j)$ for $c \in \{1,\alpha,\alpha^2\}$ and $j \in \{\alpha,\alpha^2\}$.  Since $y$ is fixed,
$\sigma(e(y,c,j)) = e(y,c^2,j^2)$. There are three orbits of size~2, so $C(x)+a \sim [2,2,2]$.

Now suppose that $\sigma$ fixes exactly one or two roots of $x^{q+1}+ax+a$,
so there is at least one rational root $r$. Select $r_0$ to be any root
that is not fixed by $\sigma$, and let $r_1 = \sigma(r_0)$.  Then
$\sigma(r_1)\not\in\{r,r_1\}$, so $\sigma(r_1)=r_{1+b}$ with $b\in\F_4^\x$.
Set $y=(r_1-r)/(r_1-r_0)$. By Theorem~\ref{qplus1Thm2},
there is a unique $\gamma \in \PGL_2(\F_q)$ such that $\sigma(r_w) = r_{\gamma(w^2)}$ and $\sigma(y)=\gamma^{-1}(y)$.
Now $\gamma(\infty)=\infty$, $\gamma(0)=1$, and $\gamma(1)=1+b$, therefore
$\gamma = \textmatrix b 1 0 1$  and $\sigma(r_w) = r_{b w^2+1}$. If $b=1$, then $\sigma$ fixes $r_\infty$, $r_\alpha$, and $r_{\alpha^2}$,
contradicting that $\sigma$ fixes exactly one or two roots.  Thus, $b=\alpha$ or $b=\alpha^2$.  Let us assume that $b=\alpha$, as the other
case is similar.  Then
$\sigma(r_0)=r_1$, $\sigma(r_1) = r_{\alpha+1}$, $\sigma(r_{\alpha+1})=r_\alpha$, and $\sigma(r_\alpha)=r_0$, thus $x^5+ax+a\sim[1,4]$.
The action on the roots of $C(x)+a$ is given by
\begin{eqnarray*}\sigma\left(e(y,c,j)\right)&=& e(\gamma^{-1}(y),c^2,j^2)=e((y+1)/b,c^2,j^2)=e(y+1,c^2/b,j^2)\\
&=& e(y,c^2/b,j^2+c^2/b+c/b^2)
\end{eqnarray*}
for $c\in\F_q^\x$ and $j\in\F_{q,1}$.
Setting $b=\alpha$, the $\sigma$-orbits are as follows: 
$$e(y,1,\alpha) \to e(y,\alpha^2,\alpha) \to e(y,1,\alpha^2) \to e(y,\alpha^2,\alpha^2) \to e(y,1,\alpha)$$
while $e(y,\alpha,\alpha)$ and $e(y,\alpha,\alpha^2)$ are fixed.
Thus, $C(x)+a\sim [1,1,4]$, as claimed.

It remains to consider the case where $x^5+ax+a$ has no rational roots. Select three roots as follows:
$r_\infty$ belongs to an orbit of odd order, $r_0 = \sigma(r_\infty)$, and $r_1=\sigma(r_0)$.
Either $\sigma(r_1)=r_\infty$ or $\sigma(r_1)=r_c$ with $c\in \{\alpha,\alpha^2\}$. Since $\sigma(r_w)=r_{\gamma \sigma(w)}$,
$\gamma$ takes $(\infty,0,1)$ to $(0,1,\infty)$ or $(0,1,c)$. In the former case, $\gamma = \textmatrix 0 1 1 1$,
and in the latter case, 
$\gamma = \textmatrix 0 1 c 1$. In the latter case, $r_{c+1}$ is fixed,
because $\sigma(r_{c+1}) = r_{\gamma(\sigma(c+1))} = r_{\gamma(c)}  = r_{1/(c^2+1)} = r_{c+1}.$  Since we are assuming there are no rational
roots, this case can be eliminated from consideration.  Thus, $\gamma=\textmatrix 0111$.
Since $\sigma(r_w)=r_{\gamma(w^2)}=r_{1/(w^2+1)}$, the $\sigma$-orbits on
$\{r_w : w \in \P^1(q)\}$ are $(r_\infty\ r_0\ r_1)(r_\alpha\ r_{\alpha+1})$.
The action on roots of $C(x)+a$ is
\begin{eqnarray*}
\sigma(e(y,c,j)) &=& e(\gamma^{-1}(y),c^2,j^2) = e((1/y)+1,c^2,j^2) = e(1/y,c^2,j^2 + c^2 + c)\\
&=& e(y,(j^2/c^2)+1+1/c,j^2+c^2+c) = e(y,cj^2+1+c^2,j^2+c^2+c),
\end{eqnarray*}
where we used $c^3=j^3=1$ since $|\F_q^\x|=3$.  There is a single $\sigma$-orbit:
$$e(y,1,\alpha)\to e(y,\alpha^2,\alpha^2) \to e(y,\alpha,\alpha^2) \to e(y,1,\alpha^2) \to e(y,\alpha,\alpha) \to e(y,\alpha^2,\alpha) \to e(y,1,\alpha).$$
So in this case, $x^5+x+1/a \sim [3,2]$ and $x(x+1)^5 + a \sim [6]$.

\section{Dihedral group} \label{dihedralSec}

Let $\LL$ denote the splitting field of $x^{q+1}+ax+a$ and let $e\in \LL$ be a root of $C(x)+a$.  As shown in Figure~\ref{splitting_diagram2},
when $\F=\F_q(a)$ with $a$ transcendental we have $[\LL:\F(e)]=2(q+1)$.  Thus, 
$\Gal(\LL/\F(e))$ is a subgroup of order $2(q+1)$  in $\PGL_2(\F_q)$, and by \cite[Chapter~XII]{Dickson}, the only such subgroup
is a dihedral group.  In this section, we give explicit formulas for this dihedral group.  Later, we will use these formulas to give
a new proof that $C(x)$ is exceptional when $n$ is odd.  

First, we discuss $\PGL_2(\F_q)$ (where $q=2^n$) in more detail. Dickson \cite[Chapter~XII]{Dickson} showed that all nontrivial elements of 
$\PGL_2(\F_q)$ have order 2, or have
order dividing $q-1$, or have order dividing $q+1$.  In fact, he enumerated these:
\begin{itemize}
\item[] $q^2-1$ elements have order~2
\item[] $(q+1)(q/2)(q-2)$ elements have order dividing~$q-1$
\item[] $q^2(q-1)/2$ elements have order dividing~$q+1$.
\end{itemize}
Including also the trivial element,  these numbers add up correctly to the full cardinality of $\PGL_2(\F_q)$:
$$q^2-1 + (q+1)(q/2)(q-2) + q^2(q-1)/2 + 1 = q(q-1)(q+1).$$
The next lemma explicitly describes these elements in a simple manner. 

\begin{lemma}  \label{lem:orderGamma}
	Let $1\ne\gamma=\textmatrix ABCD$. Then order$(\gamma)=2$ iff $A+D=0$.
	If $A+D\ne0$ and $BC=0$, then order$(\gamma)$ divides $q-1$.
	If $A+D\ne0$ and $BC\ne 0$, then on dividing through by $(BC)^{1/2}$,
	one can assume $\gamma=\textmatrix A {1/C} C D$. Writing $D=A+1/j$
	with $j\in\F_q^\x$,
	then order$(\g)$ divides $q-1$ iff $j\in\F_{q,0}$ and 
	order$(\g)$ divides $q+1$ iff $j\in\F_{q,1}$.
\end{lemma}
\begin{proof}
$$\gamma^2= \begin{pmatrix} A^2+BC & B(A+D) \\ C(A+D) & D^2+BC \end{pmatrix}$$
and this is scalar iff $A+D=0$. Thus, order$(\gamma)=2$ iff $A=D$.

In the remaining cases, $A+D\ne0$. If one normalizes so that $AD-BC=1$
(by dividing through by the scalar $(AD-BC)^{1/2}$),
then the eigenvalues are distinct and satisfy $x^2+(A+D)x+1=0$.
Since the product of the eigenvalues is 1, they are either a rational pair
$\{w,w^{-1}\}$ or a conjugate pair $\{\z,\z^q\}$ with $\z^{q+1}=1$.
In the former case, order$(\gamma)={\rm order}(w)$ divides $q-1$,
and $1/(A+D)\in\F_{q,0}$ by Lemma~\ref{Fq1Lemma}.
In the latter case, order$(\gamma)=
{\rm order}(\zeta)$ divides $q+1$, and $1/(A+D)\in\F_{q,1}$
by Lemma~\ref{Fq1Lemma}. 
Note that $(AD-BC)/(A+D)^2$ is invariant under scalar multiplication
and it equals $1/(A+D)^2$ when $AD-BC=1$.
Thus, for any matrix $\textmatrix ABCD\in\PGL_2(\F_q)$
with $A+D\ne0$ (not necessarily normalized to have determinant 1),
$$\text{${\rm order}(\gamma)$ divides $q-1$ iff $\Tr_{\F_q/\F_2}\left((AD-BC)/(A^2+D^2)\right) = 0$;}$$
$$\text{${\rm order}(\gamma)$ divides $q+1$ iff $\Tr_{\F_q/\F_2}\left((AD-BC)/(A^2+D^2)\right) = 1$.}$$

If $BC=0$, then the absolute trace of $(AD-BC)/(A+D)^2$ is always zero because
$$\frac{AD}{A^2+D^2} = \frac D{A+D} + \left(\frac D{A+D} \right)^2.$$
Thus, order$(\gamma)$ divides $q-1$. 
If $BC \ne 0$, then we renormalize so that $BC=1$. Then
$$\frac{AD+1}{A^2+D^2} =\frac D{A+D} + \left(\frac D{A+D} \right)^2 + \frac 1{(A+D)^2}.$$ 
Set $j=1/(A+D)\in\F_q^\x$. We have proved that if
\begin{equation} \gamma = \begin{pmatrix} A & 1/C  \\ C & A + 1/j \end{pmatrix} \in \PGL_2(\F_q) \label{j_form} \end{equation}
then order$(\gamma)$ divides $q+1$ iff $j\in \F_{q,1}$, and 
otherwise order$(\gamma)$ divides $q-1$.
\end{proof}
%Further, 
%every nontrivial $\gamma$ with order dividing $q+1$ has this form with
%$j\in\F_{q,1}$.
Now we can verify Dickson's counts.  Each $\gamma$ with $A+D=0$
can be written uniquely
as either $\textmatrix 1bc1$ with $bc \ne 1$, or as $\textmatrix 0bc0$ with $bc=1$, so they are in bijection with the ordered pairs $(b,c)\in\F_q \x \F_q$.
Excluding the identity, it follows that there are exactly $q^2-1$ elements of order~2, in agreement with Dickson.

By the lemma, $1\ne\gamma\in\PGL_2(\F_q)$ has order dividing
$q+1$ iff it has the form (\ref{j_form}) with $A\in\F_q$,
$C\in\F_q^\x$, and $j\in\F_{q,1}$.
Note that $j^2\det(\gamma) = j^2 A^2 + j A + j^2$. Since 
$\Tr_{\F_q/\F_2}(j)=1$, this has no rational solutions for $A$, and so
we obtain an invertible matrix for every triple 
$(A,C,j) \in \F_q \x \F_q^\x  \x \F_{q,1}$.  These are distinct as elements of $\PGL_2(\F_q)$ because
the normalization $BC=1$ uniquely determines the scalar multiple.
There are exactly $q(q-1)(q/2)$ such
triples $(A,C,j)$, which agrees with Dickson's count for the number of elements
whose order divides $q+1$.

To directly count elements of order dividing $q-1$, there are
$(q-2)(1+2(q-1))=(q-2)(2q-1)$ elements of the form
$\textmatrix ABC1$ with $BC=0$, since $A\in\F_q\setminus\F_2$ and
$(B,C) \in \{(0,0)\} \cup \left(\F_q^\x \x \{0\}\right) \cup 
\left(\{0\} \x \F_q^\x\right)$.
If $BC\ne0$, then $\gamma$ can be written uniquely in the form~(\ref{j_form})
with $j \in \F_{q,0}\setminus\{0\}$. 
There are $q-1$ choices for $C$ and $q/2-1$ choices for~$j$.
Note that $j^2\det(\gamma)=j^2A^2+jA+j^2$
has two rational roots, namely $A=b/j$ or $A=(b+1)/j$,
where $b^2+b=j^2$. (Here $b\in\F_q$ because $j\in\F_{q,0}$.)
Thus, there are $q-2$ values for $A$ that make the determinant nonzero. 
The number of $\gamma$ of order dividing $q-1$ and having
$BC=1$ is then $(q-1)(q/2-1)(q-2)$, and the total number with order 
dividing $q-1$ is 
$(q-2)(2q-1)+(q-2)(q/2-1)(q-1)=(q-2)(q^2+q)/2$,
in agreement with Dickson's count.

By a direct computation, if $A_1+A_2+1/j\ne 0$, then
\begin{equation*} \begin{pmatrix} A_1 & 1/C  \\ C & A_1 + 1/j \end{pmatrix} 
\begin{pmatrix} A_2 & 1/C  \\ C & A_2 + 1/j \end{pmatrix} 
=\begin{pmatrix} A_3 & 1/C  \\ C & A_3 + 1/j \end{pmatrix},\quad \text{where $A_3 = \frac {1+A_1A_2}{A_1+A_2+1/j}$.}
\end{equation*}
Thus, if $C,j$ are held fixed, then the elements of the form (\ref{j_form}), 
together with the identity, form a group which we denote by
${\cal C}_{j,C}$.  (Here we must exclude the matrices of determinant zero.)  We associate the identity element
with ``$A=\infty$''.  Since the determinant is nonzero for all $A\in\F_q$
when $j\in\F_{q,1}$, and for all but two values of $A$ if $j\in\F_{q,0}$,
the order of $\calC_{j,C}$ is $q+1$ if $j\in\F_{q,1}$, and $q-1$ if
$j\in\F_{q,0}$.

%It is useful to observe that
%\begin{equation} \begin{pmatrix} C&0\\ 0&1 \end{pmatrix} \begin{pmatrix} A & 1/C  \\ C & A + 1/j \end{pmatrix}
%\begin{pmatrix} 1/C&0\\ 0&1 \end{pmatrix}   = \begin{pmatrix} A & 1  \\ 1 & A + 1/j \end{pmatrix} \label{conjC}
%\end{equation}
%Often this makes it easy to reduce to the case $C=1$.  We let ${\cal C}_j $ denote the cyclic group consisting
%of matrices $M_A=\textmatrix  A 1 1 {A+1/j}$.  The value~$j$ does not change under the above conjugation. This 
%is an advantage of the normalization $BC=1$.

%Note the formula:
%\begin{equation*} \begin{pmatrix} 1&\frac 1{jC}\\ 0&1 \end{pmatrix} \begin{pmatrix} A & \frac 1 C  \\ C & A + \frac 1 j \end{pmatrix}
%\begin{pmatrix} 1&\frac 1{jC} \\ 0&1 \end{pmatrix}   = \begin{pmatrix} A & \frac 1 C  \\ C & A + \frac 1 j \end{pmatrix}^{-1}.
%\end{equation*}
%Thus, the group ${\cal D}_{j,C}$ that is generated by ${\cal C}_{j,C}$ and 
%$\textmatrix {1\ }{1/(jC)} 0 1$ has order
%$2(q-1)$ (if $j\in \F_{q,0}$) or $2(q+1)$ (if $j\in \F_{q,1}$).

From here on, we restrict attention to elements whose order divides $q+1$,
\ie, elements of the form~(\ref{j_form}) with $j\in\F_{q,1}$.
The next lemma shows how to diagonalize ${\cal C}_{j,C}$ when $j\in\F_{q,1}$.

\begin{lemma} \label{lem:D0} Let 
$\calC_0 = \left\{\textmatrix \nu 0 0 {1/\nu} : \nu \in \mu_{q+1}\right\}
\subset \PGL_2(\F_{q^2})$.  Let $C\in\F_q^\x$ and $j\in\F_{q,1}$, and
write $j=1/\ang\zeta$ with $\zeta \in \mu_{q+1}\setminus \{1\}$.
Let $M=\textmatrix C \zeta C {1/\zeta}$. Then
$${\cal C}_{j,C} = M^{-1} {\cal C}_0 M,\qquad
	\begin{pmatrix} 1 & \frac1{jC} \\ 0 & 1 \end{pmatrix} = 
		M^{-1} \begin{pmatrix} 0&1\\1&0 \end{pmatrix}  M.$$
\end{lemma}

\begin{proof} First,
\begin{equation*} M^{-1} \textmatrix 0110 M =
\begin{pmatrix} \frac1\zeta & \zeta \\ C & C \end{pmatrix}
\begin{pmatrix} C & \frac1\zeta \\ C & \zeta \end{pmatrix} =
\begin{pmatrix} C\ang\zeta & \ang{\zeta^2} \\ 0 & C\ang\zeta \end{pmatrix}
	= \begin{pmatrix} 1 & \ang{\zeta}/C \\ 0 & 1 \end{pmatrix} 
= \begin{pmatrix} 1 & \frac1{jC} \\ 0 & 1 \end{pmatrix}.
\end{equation*}
%\begin{eqnarray*} M^{-1} \textmatrix 0110 M &=& 
%\begin{pmatrix} 1/\zeta & \zeta \\ C & C \end{pmatrix}
%\begin{pmatrix} C & 1/\zeta \\ C & \zeta \end{pmatrix} =
%\begin{pmatrix} C\ang\zeta & \ang{\zeta^2} \\ 0 & C\ang\zeta \end{pmatrix}\\
	%&=& \begin{pmatrix} 1 & \ang{\zeta}/C \\ 0 & 1 \end{pmatrix} 
%= \begin{pmatrix} 1 & \frac1{jC} \\ 0 & 1 \end{pmatrix}.
%\end{eqnarray*}
Next,
	\begin{equation*} M^{-1}\textmatrix \nu 0 0 {1/\nu} M =
\begin{pmatrix} \frac 1\zeta & \zeta \\ C & C \end{pmatrix}
	\begin{pmatrix} C\nu & \nu \zeta \\ C/\nu & 1/(\zeta\nu) \end{pmatrix} 
= \begin{pmatrix} C \ang{\zeta/\nu} & \ang\nu \\
C^2\ang\nu & C\ang{\nu\zeta} \end{pmatrix} 
	= \begin{pmatrix}  \frac{\ang{\zeta/\nu}}{\ang\nu} & 1/C \\
	C & \frac{\ang{\nu\zeta}}{\ang \nu} \end{pmatrix}. 
\end{equation*}
%\begin{eqnarray*} M^{-1}\textmatrix \nu 0 0 {1/\nu} M &=&
%\begin{pmatrix} 1/\zeta & \zeta \\ C & C \end{pmatrix}
%\begin{pmatrix} C\nu & \nu \zeta \\ C/\nu & 1/(\zeta\nu) \end{pmatrix} 
%= \begin{pmatrix} C \ang{\zeta/\nu} & \ang\nu \\
%C^2\ang\nu & C\ang{\nu\zeta} \end{pmatrix} \\
%&=& \begin{pmatrix}  \ang{\zeta/\nu}/\ang\nu & 1/C \\
%C & \ang{\nu\zeta}/\ang \nu \end{pmatrix}. 
%\end{eqnarray*}
	Write this as $\textmatrix A{1/C}CD$. By (\ref{xyangProperty}),
$$A+D=\frac{ \ang{\zeta/\nu} + \ang{\zeta\nu} } {\ang\nu} 
= \frac{ \ang\z \ang \nu} {\ang \nu} = \ang\z = 1/j.$$
This shows that $M^{-1} \textmatrix \nu 0 0 {1/\nu} M \in{\cal C}_{j,C}$.
Since $\calC_0$ and $\calC_{j,C}$ both have cardinality $q+1$
	and since conjugation by $M$ is injective,
	we conclude that $M^{-1}\calC_0M=\calC_{j,C}$.
\end{proof}

The next proposition shows that the 
subgroups of order $q+1$ are cyclic and are naturally parameterized
by $\F_{q,1} \x \F_q^\x$.

\begin{proposition}   \label{prop:dihedral}
Let $q=2^n$, $C\in \F_q^\x$, and  $j\in \F_{q,1}$. For $A\in \P^1(\F_q)$, define 
$M_A \in \PGL_2(\F_q)$ by the formula
$$ M_A = \begin{pmatrix} A & 1/C \\ C & A+1/j \end{pmatrix} \qquad\text{if $A\in\F_q$, }\qquad
M_\infty = \textmatrix 1001.$$ 
Define
\begin{equation}{\cal C}_{j,C} = \{ M_A : A \in \P^1(\F_q) \},\qquad {\cal D}_{j,C} = {\cal C}_{j,C} \cup \left\{\, M_A \begin{pmatrix} 1&\frac 1 {jC}\\0&1\end{pmatrix}
: A \in \P^1(\F_q) \,\right\}.\label{DjCDef}
\end{equation}
Then ${\cal C}_{j,C}$ is a cyclic group of order $q+1$, and ${\cal D}_{j,C}$ is a dihedral group of order $2(q+1)$. 
If $1\ne\g\in \PGL_2(\F_q)$ and order$(\gamma)$ divides $q+1$ then
$\gamma\in \calC_{j,C}$ for a unique pair $(j,C)\in \F_{q,1}\x \F_q^\x$, 
and $\gamma$ has no fixed points in $\P^1(\F_q)$.  We have
\begin{equation} \text{$M_A M_B = M_K$, with $K=\frac{1+A B}{1/j+A +B}$,}
\label{kappaEqn}
\end{equation}
\begin{equation} 
\begin{pmatrix} 1&\frac 1 {jC}\\0&1 \end{pmatrix} M_A \begin{pmatrix} 1&\frac 1{jC}\\0&1 \end{pmatrix}= M_A^{-1} = M_{A+1/j}.\label{dihedralRelation}\end{equation}
\end{proposition}

\begin{proof} 
$M_A$ is invertible because $j^2\det(M_A) = (jA)^2+(jA)+j^2 \in \F_{q,1}$.
${\cal C}_{j,C}$ is cyclic of order $q+1$ by Lemma~\ref{lem:D0}.
The relation (\ref{dihedralRelation}) is straightforward to check, 
and this verifies that ${\cal D}_{j,C}$ is a dihedral group of order $2(q+1)$.

By Lemma~\ref{lem:orderGamma}, every nontrivial $\gamma$
of order dividing $q+1$ has the form (\ref{j_form}) with
$C\in\F_q^\x$ and $j\in\F_{q,1}$, thus it belongs to a group
$\calC_{j,C}$. Since the normalization $BC=1$ uniquely 
determines the scalar multiple,
the values $j$ and $C$ are uniquely determined. 

The fact that $M_A$ has no fixed points on $\P^1(\F_q)$ when $A\in \F_q$ is shown as follows. 
If $M_A(w)=w$ with $w\in \F_q$, then $(Aw+1/C)/(Cw+A+1/j)=w$, which is equivalent to $(Cjw)^2+(Cjw)+j^2=0$. But this would imply
$\Tr_{\F_q/\F_2}(j)=0$, a contradiction. Also $M_A(\infty)=A/C \ne \infty$.
Thus, $M_A$ has no fixed points in $\P^1(\F_q)$.

Next, we show that $M_A M_B = M_K$, where $K = (1+A B)/(1/j+A+B)$.
If $A = \infty$, then $K=B$ and $M_A = \textmatrix 1001$, so formula holds; 
similarly if $B=\infty$.  If $A $ and $B$ are both finite, then
$$M_A M_B = \begin{pmatrix} 1+AB & (A+B+1/j)/C \\ (A+B+1/j)C & 1+AB + \frac{A+B+1/j}j \end{pmatrix}.$$
If $A + B + 1/j=0$, then $M_AM_B=1$ as an element of $\PGL_2(\F_q)$, 
and also $K=\infty$, so the formula holds. If $A+B+1/j\ne0$, then
on dividing through by that constant we again find that
	$M_A M_B = M_K$. 
\end{proof}

Let $r$, $r_0$, and $r_1$ be distinct roots of $x^{q+1}+ax+a$ and 
$y=(r_1-r)/(r_1-r_0) \in \calY$.
%Let $r,r_0,r_1$ be distinct roots of $x^{q+1}+ax+a$ and $y=(r_1-r)/(r_1-r_0)$. 
By Theorem~\ref{KsubsetL}, the roots of $C(x)+a$ are
$\{e(y,c,j) : c \in \F_q^\x,\ j \in \F_{q,1}\}$, where
$$e(y,c,j) = (c y^2 + y + j/c)^{q+1}/(y^q+y)^2.$$
By (\ref{gammay}), $\gamma^{-1}(y)\in\calY$ for any $\g\in\PGL_2(\F_q)$,
therefore $\{e(\gamma^{-1}(y),c,j) : c \in \F_q^\x, j \in \F_{q,1}\}$ is
also a complete set of roots of $C(x)+a$. 

\begin{proposition} \label{fixedE}
For $\gamma\in \PGL_2(\F_q)$, 
	$$ e(\gamma^{-1} (y) ,c,j) = e(y,c,j)  \iff \gamma \in {\cal D}_{d,c/d},\quad \text{where $d=\sqrt j$.}$$
Here, ${\cal D}_{j,C}$ is the dihedral group of order $2(q+1)$ that is defined by (\ref{DjCDef}).
\end{proposition}

\begin{proof}  Suppose $e(y,c,d^2)=e(\gamma^{-1}(y),c,d^2)$, and we will show that
	$\gamma\in\calD_{d,c/d}$.  Denote this value by $e_0$. By Theorem~\ref{KsubsetL}, $C(e_0)+a=0$.
	Since $d\in\F_{q,1}$, it has the form $d=1/\ang\zeta$ with $\zeta \in \mu_{q+1}\setminus \{1\}$ by Lemma~\ref{Fq1Lemma}.
 
It was shown at the beginning of Section~\ref{sec:RootsOfC} that if $e$ 
is any root of $C(x)+a$, then there is $u\in \cj \F$ and distinct 
$\zeta_0,\zeta_1\in \mu_{q+1}\setminus\{1\}$ satisfying
$$1/e = \ang{u^{q+1}},\qquad r^2 = eT(e)^2 \ang u^{q-1},
\qquad r_0^2 = eT(e)^2 \ang{\zeta_0 u}^{q-1},\qquad 
r_1^2 = eT(e)^2 \ang{\zeta_1 u}^{q-1}.$$
(In these formulas, one can replace $(u,\zeta_0,\zeta_1)$ by 
$(1/u,1/\zeta_0,1/\zeta_1)$ without affecting $e$, $r$, $r_0$, and $r_1$.)
Further, we proved in Lemma~\ref{sevenFormulas} that
$$e=e(y,c_0,d_0^2),\quad \text{where
$c_0= \frac{\ang{\zeta_0/\zeta_1}}{\ang{\zeta_0} \ang{\zeta_1}}$, \quad
	$d_0=\frac 1{\ang{\zeta_0}}$}$$
and
\begin{equation*} u = M_0(y),\qquad {\it where\ }M_0 = 
	\begin{pmatrix} c_0/d_0 & 1/\zeta_0 \\ c_0/d_0 & \zeta_0 \end{pmatrix}.
\end{equation*}

Applying this to $e=e_0$, we find $e_0=e(y,c_0,d_0^2)=e(y,c,d^2)$. By
Theorem~\ref{KsubsetL}, the
roots $e(y,c,d^2)$ are distinct when $y$ is fixed and $c,d$ vary, therefore
$c_0=c$ and $d_0=d$.  Further, $d_0=d$ implies
$\ang{\z_0}=\ang{\zeta}$, and so $\zeta_0\in\{\zeta,\zeta^{-1}\}$.
Since we are free to change $(u,\zeta_0,\zeta_1)$ to $(1/u,1/\zeta_0,1/\zeta_1)$,
we can in fact assume that $\zeta_0=\zeta$. We have shown that
$e_0=e(y,c,d^2)$ and $d=1/\ang\zeta$ imply
\begin{equation} 1/e_0=\ang{u^{q+1}}, \qquad {\it where}\ 
u = M(y),\quad  M = \begin{pmatrix} c/d & 1/\zeta \\ c/d & \zeta \end{pmatrix}.
\label{uyReln} \end{equation}
Since $\gamma^{-1}(y)\in \calY$ and $e_0=e(\gamma^{-1}(y),c,d^2)$,
the same reasoning implies
\begin{equation} 1/e_0=\ang{\widetilde u^{q+1}}, \qquad {\it where}\ 
\widetilde u = M(\gamma^{-1}(y)),\quad
	M = \begin{pmatrix} c/d & 1/\zeta \\ c/d & \zeta \end{pmatrix}.
	\label{uyReln2} \end{equation}

%Recall from Theorem~\ref{qplus1Thm2} that $\gamma^{-1}(y) = (r_{\gamma(1)}-r_{\gamma(\infty)})/
%(r_{\gamma(1)}-r_{\gamma(0)})$.
%By applying the above reasoning with $\gamma^{-1}(y)$, $r_{\gamma(\infty}$,
%$r_{\gamma(0)}$, and $r_{\gamma(1)}$ in place of $y,r,r_0,r_1$,
%there are $\widetilde u$, $\widetilde \zeta$ and $\widetilde \rho$
%such that $1/e_0 = \ang{\widetilde u^{q+1}}$ and 
%\begin{equation} e = e(\gamma^{-1} (y),\widetilde c, \widetilde d^2),\qquad
%\widetilde u = \widetilde M(\gamma^{-1}(y)), \label{uyTildeReln} \end{equation}
%where $\widetilde c$ and $\widetilde d$ have the same formula as $c_0$ and~$d_0$, but with $\zeta$ and $\rho$ replaced by 
%$\widetilde \zeta$ and $\widetilde \rho$, and where 
%$$\widetilde M = 
%\begin{pmatrix} \widetilde c/\widetilde d & 1/\widetilde \zeta \\
%\widetilde c/\widetilde d & \widetilde \zeta \end{pmatrix}.$$
%As before, we are free to replace 
%$(\widetilde u, \widetilde \zeta, \widetilde \rho)$
%by $(1/\widetilde u,1/\widetilde \zeta, 1/\widetilde \rho)$.
	%Since $e_0=e(\gamma^{-1}(y),\widetilde c, \widetilde d)$ and
	%also $e_0=e(\gamma^{-1}(y),c,d)$, we have $\widetilde c=c$ and
	%$\widetilde d = d$.
%The condition $d=\widetilde d$ implies $\widetilde \zeta \in \{\zeta,\zeta^{-1}\}$.  By possibly replacing 
%$(\widetilde u, \widetilde \zeta, \widetilde \rho)$ by 
%$(1/\widetilde u,1/\widetilde \zeta, 1/\widetilde \rho)$,
%we may arrange that $\widetilde \zeta = \zeta$.  
%Then $\widetilde M= M$.

Because $\ang{u^{q+1}} = \ang{\widetilde u^{q+1}} = 1/e_0$, there is $\nu \in \mu_{q+1}$ such that
$$\text{$\widetilde u = \nu^2 u$ or $\widetilde u= \nu^2/u$.}$$
In other words, $\widetilde u = R(u)$, where $R$ belongs to the
dihedral group ${\cal D}_0$ that is generated by matrices $\textmatrix
	\nu 00{1/\nu}$ with $\nu\in\mu_{q+1}$ and by $\textmatrix 0110$.
By (\ref{uyReln}) and (\ref{uyReln2}), 
\begin{equation}\gamma^{-1} (y) = M^{-1}(\widetilde u) = M^{-1} R (u) = M^{-1} R M (y).
\label{gammaDelta} 
\end{equation}
Let $\delta^{-1}=M^{-1} R M$. By Lemma~\ref{lem:D0},
$\delta \in {\cal D}_{d,c/d}$.
Since $y \not\in \F_{q^2}$ (by Theorem~\ref{qplus1Thm}), 
an equality $\gamma^{-1} (y) = \delta^{-1} (y)$ with 
$\gamma,\delta \in\PGL_2(\F_q)$ implies $\gamma= \delta$. Thus,
	$\gamma = \delta \in {\cal D}_{d,c/d}$. 

We have shown that 
\begin{equation}\text{$e(\gamma^{-1}(y),c,d^2)=e(y,c,d^2)$ implies $\gamma\in{\cal D}_{d,c/d}$.  } \label{implies}\end{equation}
To prove the converse, we use a counting argument. Let ${\cal E}$ denote the roots of $C(x)+a$.
For $e'\in{\cal E}$, let $H_{e'} = \{ \gamma \in \PGL_2(\F_q) : e(\gamma^{-1} (y) ,c,d) = e' \}$.
Then $\PGL_2(\F_q)$ is the disjoint union of $H_{e'}$, for $e' \in {\cal E}$, and so the average size of
$H_{e'}$ is $|\PGL_2(\F_q)|/|{\cal E}| = q(q+1)(q-1)/\deg(C) = 2(q+1)$.  On the other hand, if 
$\gamma_1,\gamma_2 \in H_{e'}$, then setting $y'=\gamma_2^{-1} (y)$, 
the formula $e'=e(\gamma_1^{-1}(y),c,d^2)=e(\gamma_2^{-1}(y),c,d^2)$ implies
$$e'=e(\gamma_1^{-1} \gamma_2 (y'),c,d^2) = e(y',c,d^2).$$
Then $\gamma_1^{-1} \gamma_2 \in {\cal D}_{d,c/d}$ by~(\ref{implies}).  Thus, $|H_{e'}| \le |{\cal D}_{d,c/d}| = 2(q+1)$.
Since the average size of $H_{e'}$ is $2(q+1)$, it must be that $|H_{e'}|=2(q+1)$.
In particular, $|H_e|=2(q+1)$. By (\ref{implies}), $H_e \subset {\cal D}_{d,c/d}$, and by
comparing cardinalities, equality must hold.
\end{proof}

\section{Exceptionality of $C(x)$} \label{sec:Exceptional}

A polynomial $P(x)\in\F_r[x]$ is said to be {\it exceptional} over a finite field $\F_r$ if it induces a permutation
on ${\cal E}$ for infinitely many extension fields ${\cal E}=\F_{r^m}$.  It was proved in \cite{CM} that $C(x)$ and some 
related polynomials are exceptional over $\F_2$ when $n$ is odd. Specifically, $C(x)$ induces a permutation on $\F_{2^m}$ iff  $(2m,n)=1$.
The first polynomials in this family were found by P.\ M\" uller \cite{Muller}
with $q=8$, degree~28. M\"uller's search was motivated by some 
deep work by Fried, Guralnick and Saxl suggesting that new examples of 
exceptional polynomials
might be found in characteristic~$p=2$ or~$p=3$
having degree $(q/2)(q-1)$, where $q=p^n$ and $n$ is odd.

In this section we give a new proof that $C(x)$ is exceptional.  We emphasize 
that the next proposition is known, and only the proof is new.

\begin{proposition} If $q=2^n$, then $C(x)=xT(x)^{q+1}$ induces a permutation on $\F_{2^m}$ iff $(n,2m)=1$. 
\end{proposition}

\begin{proof}
	If $n=1$, then $C(x)=x$ and the proposition is trivially true. 
	Now assume $n>1$.
Note that $xT(x)= x + x^2 + \ldots x^{2^{n-1}}$, and its roots are precisely $\F_{2^n,0}$. Since $C(x)=xT(x)^{q+1}$, the set of roots of $C(x)$
is also $\F_{2^n,0}$.
In particular, if $\F_{2^m} \cap \F_{2^n,0}\ne \{0\}$, then $C(x)$ is not a permutation polynomial on $\F_{2^m}$. Now $\F_{2^m}\cap \F_{2^n} = \F_{2^k}$ where 
$k=(m,n)$.  Since $\F_{2^k,0}\subset \F_{2^n,0}$, for $C(x)$ to be a permutation polynomial on $\F_{2^m}$, it is necessary that
	$\F_{2^k,0} = \{0\}$, which forces $k=1$, {\it i.e.}, $(m,n)=1$.  If $n$ is even, then $T(1)=0$ so $C(0)=C(1)=0$. Thus,
another necessary condition for $C(x)$ to be a permutation polynomial on $\F_{2^m}$ is that $n$ is odd.  Together, these necessary conditions
may be written as $(2m,n)=1$.

From here on, let $\F = \F_{2^m}$, where $(2m,n)=1$. Then
$\F\cap\F_{q,0}=(\F\cap\F_q)\cap\F_{q,0}=\F_2\cap\F_{q,0}=\{0\}$,
so 0 is the unique root of $C(x)$ in $\F$.
Let $e\in\F^\x$ and $a=C(e)$. Then $a\in\F^\x$.
To prove $C$ permutes $\F$, we must show that if $e'\in\F$ and
	$C(e')=a$, then $e=e'$.  
%Assuming this condition, then $\F_{2^m}\cap \F_q=\F_2$ and $\F_{2^m}\cap \F_{q,0} = \F_2\cap\F_{q,0}=\{0\}$, so $C(x)$ has no roots in $\F_{2^m}^\x$. Thus, $C(x)$ sends $\F_{2^m}^\x$ to $\F_{2^m}^\x$.
%To prove that $C(x)$ is a permutation polynomial on $\F_{2^m}$, 
%it suffices to show that if $e,e'\in\F_{2^m}^\x$ and $C(e)=C(e')$,
%then $e=e'$. Here $a=C(e)$ is nonzero, and $e,e'$ are roots of
%$C(x)+a$.
%
%if $e \in \F^\x$ and $a=C(e)$, then $e$ is the unique root of $C(x)+a$.  Suppose that $e'$ also satisfies $a=C(e')$, and we will show that $e=e'$.

By Theorem~\ref{KsubsetL}, 
$e = e(y_0,c_0,j_0)$ for some $c_0\in \F_q^\x$ and $j_0 \in \F_{q,1}$, 
where $y_0\in{\cal Y}$.  Since $n$ is odd, $\Tr_{\F_q/\F_2}(1)=1$, 
and so $j_0 = 1 + b_0^2 + b_0$ for some $b_0\in \F_q$. 
By Lemma~\ref{galCLemma}, $e(y_0,c_0,j_0)=e(c_0y_0,1,j_0)=e(c_0y_0+b_0,1,1)= e(y,1,1)$, where $y=c_0y_0+b_0\in {\cal Y}$. 
	By Theorem~\ref{KsubsetL}, if $e'\in\F^\x$ is another root of $C(x)+a$,
	then $e'=e(y,c,j)$ for some $c\in\F_q^\x$
and $j \in \F_{q,1}$.  Writing
$j=1+b+b^2$, we see that $e'=e(cy+b,1,1)$.  Thus,
$$e = e(y,1,1),  \qquad e' = e(cy+b,1,1).$$

Since $e$ and $e'$ are rational, they are fixed by every element of $\Gal(\LL/\F)$, where $\LL=\F \circ \F_q(y)$ is the splitting field of $C(x)+a$.
Since $\LL$ is finite, $\Gal(\LL/\F)$ is generated by the Frobenius element:
	$$\sigma(u) = u^{|\F|}=u^{2^m}.$$
If $w \in \F_q$ and $\sigma(w)=w$, then 
	$w\in \F\cap\F_q=\F_2$. Further, $\sigma(w)+w=1$ is impossible, since
	$\Tr_{\F_q/\F_2}(\sigma(w)+w)=0$ and $\Tr_{\F_q/\F_2}(1)=1$. 
	Together, these observations imply:
\begin{equation} \text{\it If $w\in\F_q$ and $\sigma(w)+w\in\F_2$, then $w\in\F_2$.}
	\label{winF2}
\end{equation}

Let $\gamma$ be the element of $\PGL_2(\F_q)$ such that $\sigma(y)=\gamma^{-1}(y)$. (Such $\gamma$ exists and is unique by Theorem~\ref{qplus1Thm2}.)
	Then $\sigma(e) = e(\gamma^{-1}(y),1,1)$, 
	$\sigma(e') = e(\sigma(c) \gamma^{-1}(y) + \sigma(b),1,1)$.
Since $\sigma$ fixes $e$ and $e'$, 
$$e=e(y,1,1) = e(\gamma^{-1}(y),1,1), 
\qquad e'=e(cy+b,1,1) = e(\sigma(c)\gamma^{-1}(y)+\sigma(b),1,1).$$
By Proposition~\ref{fixedE}, the first equality implies that 
$\gamma \in {\cal D}_{1,1}$, and the second implies 
	$$\sigma(c) \gamma^{-1}(y) + \sigma(b) = \delta^{-1}(cy+b),$$
where $\delta \in {\cal D}_{1,1}$.  Equivalently,
$$\begin{pmatrix} \sigma(c) & \sigma(b) \\ 0 & 1 \end{pmatrix}
\gamma^{-1}(y) = \delta^{-1} 
\begin{pmatrix} c&b\\0&1\end{pmatrix} (y).$$
Since $y \not\in\F_{q^2}$ by Theorem~\ref{qplus1Thm}, an equality
$\gamma_1(y)=\gamma_2(y)$ with $\gamma_1,\gamma_2 \in \PGL_2(\F_q)$ implies 
$\gamma_1=\gamma_2$.  Thus,
\begin{equation} \delta^{-1} = 
	\begin{pmatrix} \sigma(c) & \sigma(b) \\ 0&1 \end{pmatrix} \gamma^{-1}
	\begin{pmatrix} c & b \\ 0&1 \end{pmatrix}^{-1}. 
		\label{MAMB}
\end{equation}
			
	Since $\g,\delta \in {\cal D}_{1,1}$, we may write 
$$\gamma^{-1} = M_A \begin{pmatrix} 1& \varepsilon_1\\ 0&1\end{pmatrix},\qquad \delta^{-1} = M_B \begin{pmatrix} 1& \varepsilon_2\\ 0& 1\end{pmatrix}$$
where $\varepsilon_1,\varepsilon_2\in\F_2$ and
$$M_A = \begin{pmatrix} A & 1 \\ 1 & A+1 \end{pmatrix}\qquad \text{if $A\in \F_q$, } M_\infty = \begin{pmatrix} 1&0 \\ 0 & 1 \end{pmatrix}.$$
Note that the bottom left corner of $\gamma^{-1}$ is zero iff $A=\infty$.

First, if $A = \infty$ (so $M_A$ is the identity), 
then the relation (\ref{MAMB}) is
$$M_B = \begin{pmatrix} \sigma(c) & \sigma(b) \\ 0 & 1 \end{pmatrix} \begin{pmatrix} 1&\varepsilon_1 \\ 0 & 1 \end{pmatrix} {\textmatrix cb01}^{-1}\begin{pmatrix} 1&\varepsilon_2 \\ 0 & 1 \end{pmatrix}.$$
Since the bottom left entry is 0, $M_B = \textmatrix 1001$, and the 
relation simplifies to
$$\begin{pmatrix} 1&0 \\ 0&1 \end{pmatrix} =
\begin{pmatrix} \sigma(c) & 
\sigma(c)(\varepsilon_2+b) + c\left(\varepsilon_1 \sigma(c) 
+ \sigma(b)\right) \\ 0 & c \end{pmatrix}.$$
Since $c\in\F_q^\x$ and $\sigma(c)=c$, (\ref{winF2}) implies $c=1$.
Setting $c=1$, the top-right entry is $0=\varepsilon_1 + \varepsilon_2 + b + \sigma(b)$, and so $b + \sigma(b) \in \F_2$.  
By~(\ref{winF2}), $b\in\F_2$.  Then 
$e'=e(y+b,1,1)=e(y,1,1+b^2+b)=e(y,1,1)=e$, as we wanted to show.

Now suppose $A\in \F_q$.  Then 
$$M_B = \begin{pmatrix} \sigma(c) & \sigma(b) \\ 0 & 1 \end{pmatrix} M_A \begin{pmatrix} 1&\varepsilon_1 \\ 0 & 1 \end{pmatrix} {\textmatrix cb01}^{-1}\begin{pmatrix} 1&\varepsilon_2 \\ 0 & 1 \end{pmatrix}.$$
The right side is
$$\begin{pmatrix} \sigma(c) A + \sigma(b)\  & \left(\sigma(c) A + \sigma(b)\right)(\varepsilon_1 c + b + \varepsilon_2) + \left(\sigma(c)+\sigma(b) A + \sigma(b)\right)c \\
1 & \varepsilon_1 c + b + \varepsilon_2 + Ac + c \end{pmatrix}.$$
To be of the form $M_B$, we require that the top-right entry is 1 and the trace is 1.  This gives the two equalities:
\begin{equation} \left(\sigma(c) A + \sigma(b)\right)(\varepsilon_1 c + b + \varepsilon_2) + \left(\sigma(c)+\sigma(b) A + \sigma(b) \right)c  = 1 \label{star} \end{equation}
\begin{equation} (c + \sigma(c)) A + (\varepsilon_1+1) c + \sigma(b) + b + \varepsilon_2 = 1. \label{star2} \end{equation}

	First, if $c=1$, then (\ref{star2}) implies that $\s(b)+b\in\F_2$.
	By (\ref{winF2}), $b\in\F_2$. Then $e'=e(y+b,1,1)=e(y,1,1+b^2+b)=e(y,1,1)$ as required.

Next assume $c\ne1$, and we will obtain a contradiction. By (\ref{winF2}),
$\sigma(c)+c\not \in\F_2$, so in particular $\sigma(c)+c\ne0$.
Solve for $A$ in (\ref{star2}) and use this to eliminate $A$ in (\ref{star}).
After simplifying, we obtain:
$$\sigma(b)^2 c + b^2 \sigma(c) + \sigma(b) c + b \sigma(c) + (b + \varepsilon_1 + \sigma(b) + \varepsilon_2) c \sigma(c) +
(c+\sigma(c)) (c \sigma(c)+1) = 0.$$
On dividing  by $c \sigma(c)$, we find that $x + \sigma(x) = \varepsilon_1 + \varepsilon_2 \in \F_2$, where
$$x = b^2/c + b/c + b + c + 1/c.$$
By~(\ref{winF2}), $x\in\F_2$.

If $x=0$, then $b^2 + b(c+1) + c^2 + 1 = 0$.  Since we are assuming $c\ne 1$, 
we may divide through by $(c+1)^2$ and this gives
$$ \left( \frac b {c+1} \right)^2 + \frac b {c+1} + 1 = 0,$$
contradicting that $\Tr_{\F_q/\F_2}(1)=1$.

If $x=1$, then $b^2 + bc+b + c^2 + 1 = c$.  If $b=c$, then $c=1$, contradicting
the hypothesis that $c\ne 1$.
So we may divide through by $(b+c)^2$, and we find
$$0 = \frac {b^2 + c^2 + b + c + 1 + bc } {(b+c)^2} = 1 + \frac {c+1}{b+c} + \left(\frac {c+1}{b+c}\right)^2.$$
On taking the trace, we obtain $0=\Tr_{\F_q/\F_2}(1)$, a contradiction since $n$ is odd.
The contradiction shows that $c=1$, and we already showed that implies $e=e'$.  We conclude that 
$C(x)$ is indeed a permutation polynomial when $(n,2m)=1$.
\end{proof}

\end{document}